\documentclass{behrangstyle}
\usepackage{fancyhdr}

\usepackage[dvips]{graphicx}
\usepackage{latexsym}
\usepackage{amsmath}
\usepackage{amssymb}
\usepackage{here}
\usepackage{color}
\usepackage{subfigure}
\usepackage{amsthm}
\usepackage{multirow}
\usepackage{rotating}
\theoremstyle{definition}
\newtheorem{definition}{Definition}
\newtheorem{theorem}{Theorem}
\newtheorem{lemma}{Lemma}

\newtheorem{corollary}{Corollary}

\newtheorem*{example*}{Example}

\def\proof{\par\noindent{\scshape \underline{Proof:}} \\}
\def\done{\hfill$\blacksquare$}

\begin{document}

\pagestyle{empty}

\pagestyle{empty} \vspace{20cm}
\begin{center}
\Large{}
\end{center}
\begin{center}
\LARGE{ \bf Exploring connectivity of random\\ subgraphs of a
graph.}\\
\vspace{1mm}

 \normalsize Connectivity of random subgraphs of Cartesian products of $K_2$, $K_3$, and $P_3$. A survey of uniformly most reliable networks.
\end{center}

\vspace{1mm}
\begin{center}
\end{center}

\vspace{2mm}
\begin{center}
\large{
BEHRANG MAHJANI\\
Supervisor: Jeffrey Steif}
\end{center}

\vspace{140mm}
\begin{center}
\large{
Department of Mathematical Sciences\\
CHALMERS UNIVERSITY OF TECHNOLOGY\\
Division of Engineering Mathematics\\
G\"{o}teborg, Sweden, 2010}
\end{center}
\newpage
\vspace{5cm} {\noindent Exploring connectivity of random subgraphs
of a
graph.\\
\noindent Connectivity of random subgraphs of Cartesian products of $K_2$, $K_3$, and $P_3$. A survey of uniformly most reliable networks.\\
\noindent BEHRANG MAHJANI \vspace{5mm}\\
\noindent \copyright BEHRANG MAHJANI, 2010. \vspace{5mm}\\
\noindent Department of Mathematical Sciences\\
\noindent Chalmers University of Technology\\
\noindent SE-412 96 G\"{o}teborg\\
\noindent Sweden\\
\noindent Telephone + 46 (0)31-772 1000\\
\newpage
{\noindent Exploring connectivity of random subgraphs of a
graph\\
\noindent Connectivity of random subgraphs of Cartesian products of $K_2$, $K_3$, and $P_3$. A survey of uniformly most reliable networks.\\
\noindent BEHRANG MAHJANI\\
\noindent Department of Mathematical Sciences\\
\noindent Chalmers University of Technology \vspace{1.5cm}

\noindent \textbf{Abstract}\\\\ This work is divided into two main
parts. The first part is devoted to exploring the connectivity of
random subgraphs of cartesian products of $K_1$, $K_2$, and $P_3$.
In the second part, the author presents a short review of the
results about network
reliability.\\

The cartesian product of $K_2$, the complete graph with 2
vertices, is the cube graph $Q^n$. A random subgraph of $Q^n$,
$Q^n_{p_n}$, contains all vertices of $Q^n$, and each edge of
$Q^n$ independently with probability $p_n$. One can call $p_n$ the
percolation parameter. The author explains in detail that for
$p_n\geq1-(1/2)(\log n)^{1/n}$, $Q^n_{p_n}$ has no components with
size larger than $1$ and smaller than $2^n$, as $n\rightarrow
\infty$. It is also explained that for
$p_n=1-(1/2)\lambda^{1/n}(1+o(1/n))$, the probability that
$Q_{p_n}^n$ has no isolated point, as $n\rightarrow \infty$, tends
to $e^{-\lambda}$; hence, the probability that $Q_{p_n}^n$ is
connected tends to $e^{-\lambda}$. For constant percolation values
larger than $1/2$, when $n$ tends to infinity, almost every random
subgraph of $Q^n$ is connected; for percolation values smaller
than $1/2$, when $n$ tends to infinity, almost no random subgraph
of $Q^n$ is connected; and for percolation values equal $1/2$,
when $n$ tends to infinity, the probability that $Q^n_{1/2}$ is
connected tends to $e^{-1}$. At the end of this section, a
comparison between connectivity of a typical random graph with $M$
edges and $N$ vertices, $G\in G(N=2^n,M=n2^{n-1})$, and $Q^n$, after percolation with the parameter $p_n$ is presented.\\

This work continues with exploring the threshold function for the
cartesian product of $K_3$, the complete graph with 3 vertices,
denoted by $^3Q^n$. It is shown that for $p_n\geq
1-(1/\sqrt3)(\log n)^{1/n}$, $^3Q_{p_n}^n$ has no components with
size larger than $1$ and smaller than $3^n$, as $n\rightarrow
\infty$. Then, it is proved that for
$p_n=1-(1/\sqrt3)\lambda^{1/2n}(1+o(1/n))$, the probability that
$^3Q_{p_n}^n$ has no isolated point, as $n\rightarrow \infty$,
tends to $e^{-\lambda}$; hence, the probability that $Q_{p_n}^n$
is connected tends to $e^{-\lambda}$. At last, the author suggests
that the threshold value for connectivity of the cartesian product
of $P_3$, where $P_3$ is a path with length $2$, denoted by
$P_3^n$, is $2-\sqrt2$. One can show that for percolation values
smaller than $2-\sqrt2$, almost no random subgraph of $P_3^n$ is
connected, and for percolation values larger than $2-\sqrt2$,
almost every random subgraph of $P_3^n$ has no isolated point. The
author also shows that for percolation values larger than $0.68$
almost all random
subgraphs of $P_3^n$ are connected.\\

The last part of this work, sheds light on reliability of
networks. The main question in this part is: one is given 2
parameters, $n$ and $m$ where $n$ and $m$ are positive integers.
Among all graphs with $n$ vertices and $m$ edges, which graph $G$,
if any, maximizes the probability that when one does percolation
on $G$ with the parameter $p_n$, for all $p_n$ in $(0,1)$ there is
one component? $G$ would be called the uniformly optimally
reliable graph (UOR graph) for the parameter $n$ and $m$. It is
shown in this part, for some $m$ and $n$ there is no UOR graph,
since the graph which maximizes the probability of connectivity
depends on $p_n$ in that family of graphs. A review of results
about when the UOR graph exists is presented in this part.\\\\
\vspace{5cm}\noindent \textbf{Keywords:} Random subgraphs,
percolation, n-cube, path graph, reliable networks.
\newpage
\begin{center}
{\large {\bf ACKNOWLEDGMENTS}}\end{center}

I would like to express my deepest gratitude to Professor Jeffrey
Steif for his supervision and guidance in this thesis. His
invaluable comments and suggestions were of enormous help during
this research work. My special thanks are due to my parents who
supported me and made it possible for me to continue my studies.

\vspace{5mm}

\tableofcontents

\setcounter{page}{0} \pagestyle{fancy} \setlength{\parskip}{1ex
plus 0.5ex minus 0.2ex} \lhead{}
\chapter{Introduction}
Exploring the connectivity of random subgraphs of different
families of graphs is one of the most interesting topics in random
graphs and percolation theory. A random subgraph of a graph
$G(V_n,E_m)$ is a graph which contains all vertices of $G$, and
each edge of $G$ independently with probability $p_n$. $p_n$ is
known as the percolation parameter. Connectivity of a random
subgraph of a graph can be investigated both for small and
considerably large (when $n$ tends to infinity) graphs. For
considerably large graphs, the first step in exploring the
connectivity is to calculate $p_c$ which for all constant
$p=p_n,p\in(0,1)$ and $p<p_c$, as $n$ tends to infinity, almost
all random subgraphs of $G(V_n,E_m)$ is connected; but for all
$p\in(0,1)$ and $p>p_c$, as $n$ tends to infinity, almost no
random subgraphs of $G(V_n,E_m)$ is connected. The second step is
to investigate what happens when $p=p_c$. A more complete approach
is to calculate $p_c$ when $p_c$ depends on $n$.

For small graphs, it is of interest to find the uniformly
optimally reliable graph (UOR graph). Consider $G(n,m)$ as the
family of graphs with $n$ vertices and $m$ edges. The UOR graph is
the graph $G\in G(n,m)$ that maximizes the probability that $G$ is
connected after percolation with the parameter $p_n$ for fixed
$n,m$ and all $p_n\in(0,1)$.

One of the interesting graphs for analyzing the connectivity of
its random subgraphs is the cube graph. The cube graph $Q^n$, is a
graph with the vertices labeling $1,2,3,4,...,2^{n}-1$. Two
vertices in this graph are adjacent if their binary representation
differs only in one digit. Another way to define $Q^n$ is using
the cartesian products of $n$ copies of $K_2$, where $K_n$ is the
complete graph with $n$ vertices. It is shown by Paul Erd\"{o}s
and Joel Spencer \cite{evol} that $p_c=1/2$ for $Q^n$. An
extension of $Q^n$ is the graph with the vertices labeling
$1,2,3,4,...,3^{n}-1$, where two vertices in this graph are
adjacent if their ternary representation differs only in one
digit. We call this graph 3-cube denoted by $^3Q^n$. One can show
that $^3Q^n$ is the cartesian product of $K_3$. Connectivity of
random subgraphs of the cartesian product of $K_i$ is investigated
by Lane Clark \cite{product}. Another extension of $Q^n$ is the
cartesian product of $n$ copies of a path with length 2 which we
call it $P_3^n$.

This thesis is divided into two main parts. The first part
(chapter 4,5,6) is devoted to the connectivity of random subgraphs
of some considerably large graphs, and the second part is the
connectivity of random subgraphs of small graphs (chapter 7).
Chapter 2 presents a very short review of definition and results
in graph theory. Chapter 3 is a short review of the definitions in
random graph theory; it is explained briefly in this chapter that
how small components construct a giant component and gradually a
graph becomes connected by adding more edges to it. In chapter
$4$, the results by Bela Bollob\'{a}s \cite{bol} on finding $p_c$
for connectivity of random subgraphs of $Q^n$ is explained in
detail. In chapter $5$ the author calculates $p_c$ for
connectivity of random subgraphs of $^3Q^n$. After calculating the
threshold value for the connectivity of random subgraphs of
$^3Q^n$, the author found that this problem is solved for a
general case of the random subgraphs of cartesian product of $K_i$
\cite{product}. Chapter $6$ is an approach to find $p_c$ for
connectivity of random subgraphs of $P_3^n$. This work finishes
with chapter $7$ which is a review of the results on finding the
UOR graph. In this chapter it is shown by the author that for some
$m$ and $n$ there are no UOR graph.

\chapter{Graph theory background}
This chapter presents s short review of basic definitions in graph
theory. Most of the definitions in this chapter are extracted from
\cite{Graph}.
\section{Graph models and their matrix representation}

\subsection*{Graphs}
\begin{definition} \textbf{Graph:} A \emph{graph}
$G(V,E)$ is an ordered pair consisting of the set of
\emph{vertices} $V$, and the set of \emph{edges} $E$. Each edge is
associated with a set of vertices which are called
\emph{endpoints}. Two vertices are \emph{adjacent} if they are
joined by an edge. Two edges are \emph{adjacent} if they have a
common endpoint. A vertex is \emph{incident} to an edge and
viceversa, if that vertex is an endpoint of the edge. A
\emph{self-loop} is an edge which joins a vertex to itself. A
\emph{multi-edge } is a set of two or more edges having the same
endpoints. A \emph{simple graph} is a graph without self-loops and
multi-edges.
\end{definition}

\subsection*{Degrees}
After defining a graph, it is of interest to get familiar with the
characteristics of different graphs in order to compare them. One
of the basic characteristics is the degree of each vertices.
\begin{definition}\textbf{Degree:} The \emph{degree} of a vertex, denoted by $deg(v)$, is the
number of edges incident on that vertex plus two times the number
of its self-loops. The \emph{smallest degree} in a graph is
denoted by $\delta_{min}$ or $\delta$, and the \emph{largest
degree} in a graph is denoted by $\delta_{max}$ or $\Delta$. The
\emph{degree sequence} of a graph is the non-increasing sequence
of vertice degrees.
\end{definition}

The first question that comes into mind, after defining the degree
sequence of a graph, is if there exists a degree sequence of a
graph for each sequence of positive integers.
\begin{definition} \textbf{Graphic:} A sequence of positive integers is \emph{graphic} if there is
a permutation of it that is the degree sequence of a simple graph.
An explicit sufficient and necessary condition for a sequence of
positive integers to be graphical is:\end{definition}
\begin{theorem}
A sequence of non-negative integers $(d_1,d_2,...,d_n)$ is
graphical if and only if
\begin{align}
\sum_{i=1}^k deg(i) \leq k(k-1)+\sum_{j=k+1}^n min(k,deg(i))
\end{align}
for each $1\leq k \leq n$ \cite{Syn}.
\end{theorem}

\subsection*{Graph models}
There are many types of graphs. Some of the most important types
of simple graphs are:
\begin{definition} \textbf{Common families of graphs:} A \emph{complete graph} $K_n$ is a simple graph with $n$ vertices which every pair
of vertices is connected by an edge. A \emph{bipartite graph} $G$
is a graph with the set of vertices that can be partitioned into
two subsets $U$ and $W$, such that each edge in $G$ has one
endpoint in $U$ and one endpoint in $W$. A \emph{regular graph} is
a graph where all vertices have the same degree. A path graph is a
simple graph with $|V|=|E|+1$ that can be drawn, such that all
vertices and edges are in a single straight line. A \emph{path
graph} with $|V|=n$ and $|E|=n-1$ is denoted by $P_n$. A
\emph{hypercube graph (cube graph)} is a simple n-regular graph
with the set of vertices labels from $0$ to $2^n-1$, in which two
vertices are adjacent if their binary representation differs only
in one digit.
\end{definition}
\begin{definition} \textbf{Subgraphs:} $H$ \emph{subgraph} of $G$ is
a graph whose vertices and edges are in $G$. If $V_G=V_H$ then the
subgraph $H$ is said to \emph{span} the graph $G$. The
\emph{induced subgraph} on $U\subseteq V_G$ of $G$ is the graph
whose set of vertices is $U$, and set of edges is all edges of $G$
with two endpoints in $U$. The \emph{induced subgraph} on
$D\subseteq E_G$ of $G$ is the graph whose its set of edges is
$D$, and its set of vertices is all vertices which are incident
with an edge in $D$. A maximal connected subgraph of a graph $G$
is a \emph{component} of $G$.
\end{definition}
\begin{definition} \textbf{Cartesian product of a graph:} $G\times
H$, the \emph{Cartesian product} of $G$ and $H$ is the graph with
the set of vertices $V_G\times V_H$ and the set of edges
$(V_G\times E_H)\cup (E_G\times V_H)$.
\end{definition}

Defining a walk on a graph can help us to define some important
characteristic of a graph such as connectivity of a graph and
spanning trees.
\begin{definition} \textbf{Walk:} In a graph $G$, a \emph{walk} from vertex $v_0$ to vertex
$v_n$ is an ordered sequence
\begin{align}
W=<v_0,e_1,v_1,e_2,...,v_{n-1},e_n,v_n>
\end{align}
of vertices and edges, such that the endpoints of $e_i$ is
$\{v_{i-1},v_i\}$ for $i=1,...,n$. For a simple graph one can
abbreviate the representation as a vertex sequence
\begin{align}
W=<v_0,v_1,...,v_n>
\end{align}
\end{definition}

\begin{definition} \textbf{Tree, spanning tree:} A \emph{path} is a walk
with no repeated vertices (except the initial and final vertices).
A \emph{cycle} is a nontrivial closed path. A \emph{tree} is a
connected graph without cycle. A spanning tree of a graph is a
subgraph of a graph which is a tree.
\end{definition}

\subsection*{Matrix representations}
The last important concept in this section is that each graph can
be presented as a matrix, as follows:
\begin{definition} \textbf{Matrix representation of a graph:} The
\emph{adjacency matrix} of a simple graph $G$ is:
\begin{math}
A_G[u,v]=\left\{%
\begin{array}{ll}
    1, & \hbox{if $u$ and $v$ are adjacent;} \\
    0, & \hbox{otherwise.} \\
\end{array}%
\right.
\end{math} for all pairs of vertices $u$ and $v$ in $V_G$.
\end{definition}

\section{Connectivity}
Connectivity of a graph is one of the most important property of a
graph. A graph is \emph{connected} if for every pair of vertices
$u$ and $v$, there is walk from $u$ to $v$. There are different
types of connectivity for a graph:
\begin{definition} \textbf{Vertex-connectivity:} $\kappa_v(G)$,
\emph{vertex connectivity} of a connected graph $G$, is the
minimum number of vertices which its removal will disconnect $G$
or reduce it to a single vertex graph. A graph $G$ is
\emph{k-connected} if $\kappa_v(G)\geq k$.
\end{definition}

\begin{definition} \textbf{Edge-connectivity:}$\kappa_e(G)$,
\emph{edge connectivity} of a connected graph $G$, is the minimum
number of edges which its removal will disconnect $G$. A graph $G$
is \emph{k-edge-connected} if $\kappa_e(G)\geq k$
\end{definition}

\begin{definition}\textbf{Algebraic connectivity:} The \emph{Laplacian matrix} of
a graph is $L:=D-A$, where $A$ is the adjacency matrix and the $D$
is the diagonal matrix of vertex outdegrees. The algebraic
connectivity of an undirected graph with the Laplacian matrix $L$
is the second smallest eigenvalue of $L$. If one arranges the
eigenvalues of $L$ as : $\lambda_1(L) \leq \lambda_2(L)
\leq...\lambda_n(L) $, then $\lambda_2(L)$ is the algebraic
connectivity of a graph. The following theorem presents some
applications of algebraic connectivity:
\end{definition}
\begin{theorem}
For an undirected graph with minimum vertex degree $\delta$ and
maximum vertex degree $\Delta$, we have:
\begin{itemize}
    \item $\lambda_2 \geq 0$ with the inequality strict if and only if the
    graph is connected.
    \item $\lambda_2 \leq \frac{n}{n-1}\delta \leq \frac{n}{n-1} \triangle \leq
    \lambda_n$.
\end{itemize}
\end{theorem}

\chapter{Random graph theory background}
A \emph{random graph} is a graph with a specific number of
vertices which adjacency between two vertices are determined in a
random way \cite{Graph}. In this chapter, we define Erd\~{o}s
R\'{e}nyi random graph, and then we explain briefly some
properties of it.
\section{Evolution of Erd\~{o}s R\'{e}nyi graphs}
$G(n,M),0\leq M\leq \binom{n}{2}$, is the equiprobable space of
all simple graphs with the vertex set $V=\{1,2,...,n\}$ and $M$
edges. $G(n,p), 0< p < 1$, is the collection of all graphs with
the vertex set $V=\{1,2,...,n\}$ in which two vertices are
connected independently with the probability $p$ \cite{Graph},
\cite{bol}.

It is of interest to study the global structure of a random graph
of order $n$ (with $n$ vertices) and size M(n) (with $M(n)$
edges). Let us define $L_j(G)$ as the order of the jth largest
component of a graph $G$, where if $G$ has fewer than j components
then $L_j(G)=0$. Consider the random graph process
$\widetilde{G}=(G_t)_{t=0}^N$ where $G_t$ is getting larger by
adding more and more edges. When  $t \sim \frac{1}{2}cn$ and $c<1$
then in a.e $G_t$ the maximum of the order of its components is of
order $\log n$. When $c=1$, in a.e $G_t$ $L_1(G_{\lfloor
n/2\rfloor})$ has order $n^{2/3}$. When $t$ passes $n/2$, $L_1(G)$
begins to grow suddenly and the giant component, which is a
component whose order is much larger than other components,
appears. Eventually, small components join the giant component and
the graph becomes connected. Erd\~{o}s R\'{e}nyi proved that
$(n/2)\log n$ is the sharp threshold for connectedness \cite{bol}.
The following theorem illustrates this fundamental result:
\begin{theorem} Let $c\in \mathbb{R}$ be fixed and let $M=(n/2)\{\log n+c+o(1)\}\in
\mathbb{N}$ and $p=\{\log n+c+o(1)\}/n$. Then \cite{bol}:
\begin{align}
\textbf{P}(G_M \text{ is connected})\rightarrow e^{-e^{-c}} \text{
as } n\rightarrow \infty
\end{align}
and
\begin{align}
\textbf{P}(G_p \text{ is connected})\rightarrow e^{-e^{-c}} \text{
as } n\rightarrow \infty.
\end{align}

\end{theorem}
\noindent For more information regarding random graphs one can
check \cite{bol}, \cite{bol2}, \cite{kolchin}, \cite{janson}.
\section{Properties of almost every graphs}
A graph property $T$ is true for almost every (all) graph if for
fixed $p=p(n)$ \cite{graph}:
\begin{align}
\lim_{n\rightarrow \infty}\textbf{P}(G \in G(n,p) \text{ and } G
\text{ has property T})=1 \text{ for } p>0
\end{align}
Some of the important "almost every graph properties" are:
\begin{theorem}
For any integer $r\geq 1$ and all $p\in(0,1)$, almost every graph
contains $K_r$ \cite{Graph}.
\end{theorem}
\begin{theorem}
Almost every graph is connected for all $p\in(0,1)$ \cite{Graph}.
\end{theorem}
\begin{theorem}
For $k\in\mathbb{N}$ and all $p\in(0,1)$, almost every graph is
k-connected \cite{Graph}.
\end{theorem}

\section{Probabilistic methods}
Usually, the goal in the probabilistic method is to prove the
existence of a combinatorial structure with a certain property.
The usual approach in these methods is to first construct a
suitable probability space, then show that there exists a random
object in that space with the desired properties \cite{probab}.
Sometimes it is not easy to find the desired object, instead one
proves that there is an object which \emph{almost} satisfies the
desired conditions \cite{probab2}. Usually, it is possible to
modify the almost close object in a deterministic way so that one
gets the desired object. Markov's inequality, and Chebyshev
inequality are two important inequalities used for this purpose.
An important concept in probabilistic methods is the definition of
threshold function, which is:
\begin{definition}
$r(n)$ is called a threshold function for a graph property $T$ for
$G(n,M(n))$
if:\\
 1. When $\lim_{n\rightarrow \infty} \frac{M(n)}{r(n)}=0$  almost every graphs do not
satisfy $T$. \\2. When $\lim_{n\rightarrow \infty}
\frac{M(n)}{r(n)}=1$  almost every graphs satisfy
$T$.\\
\end{definition}

\chapter{Connected random subgraphs of the cube}
\footnote{The proof presented in this chapter
 is based on the proof presented by B.Bollob\'{a}s in \cite{bol} p.384-393.}

The cube graph, $Q^n$, is a graph with $2^n$ vertices. If one
labels each vertex of $Q^n$ from $0$ to $2^n-1$, then two vertices
are adjacent if their binary representation differs only in one
digit. Hence, one can conclude that each vertex in $Q^n$ is
connected to $n$ other vertices. In other words, $Q^n$ has
$n2^{n-1}$ edges. A random subgraph of $Q^n$ is denoted by
$Q_{p_n}^n$. $Q_{p_n}^n$ contains all vertices of $Q^n$, and each edge of
$Q^n$ independently with probability $p_n$.

It is of interest in this chapter to explore a critical value
$p_c$, which for fixed values of $p$ if $p<p_c$ then the
probability that $Q_{p_n}^n$ is connected, as $n \rightarrow \infty$,
tends to $0$; but if $p>p_c$ then the probability that $Q_{p_n}^n$ is
connected, as $n \rightarrow \infty$, tends to $1$. Burtin proved
that this critical value is $1/2$ \cite{Burtin}. Later,
P.Erd\"{o}s and J.Spencer proved that for $p=1/2$ the probability
that $Q_{p_n}^n$ is connected, as $n \rightarrow \infty$, tends to
$e^{-1}$ \cite{evol}.

 In the first section of this chapter, first the probability that
 $Q_{p_n}^n$ has no isolated point as $n\rightarrow \infty$, for fixed $p$, is investigated. It is proved that for $p<1/2$
 the probability that $Q_{p_n}^n$ has no isolated point, as $n \rightarrow \infty$, tends to
 $0$. Therefore, for $p<1/2$ the probability that $Q_{p_n}^n$ is connected, as $n \rightarrow \infty$, tends to
 $0$. Then it is proved that, for $p>1/2$ the probability that $Q_{p_n}^n$ has no isolated
point, as $n\rightarrow \infty$,  tends to $1$. In the next step,
the probability that $Q_{p_n}^n$ has no isolated point, as
$n\rightarrow \infty$, when $p$ depends on $n$ and it is close
$1/2$, is explored. It is proved that for $\lambda(n)=\lambda>0$
and $p_n=1-(1/2)\lambda^{1/n}(1+o(1/n))$, the probability that
$Q_{p_n}^n$ has no isolated point, as $n\rightarrow \infty$, tends to
$e^{-\lambda}$ \cite{bol}. Finally, it is proved that, for fixed
$p=1/2$ the probability that $Q_{p_n}^n$ has no isolated point, as
$n\rightarrow \infty$, tends to $e^{-1}$. These results are based
on P.Erd\"{o}s and J.Spencer's work \cite{evol}.

In the second section, one sheds light on the Isoperimetric
problem, which is the problem of finding an inequality which
relates the size of a subgraph to the size of its boundary. The
solution to this problem for $Q_{p_n}^n$ is presented by S.Hart
\cite{hart}. One needs such an inequality to explore the
probability that $Q_{p_n}^n$ has a component which is not the whole graph.

In the last section, the Isoperimetric inequality is applied to
prove that when $p$ depends on $n$ and $p_n\geq1-(1/2)(\log
n)^{1/n}$, then the probability that there are no components with
size larger than $1$ and smaller than $2^n$ in $Q_{p_n}^n$ , as
$n\rightarrow \infty$, tends to $1$. Therefore, for
$p_n=1-(1/2)\lambda^{1/n}(1+o(1/n))$, the probability that $Q_{p_n}^n$
is connected, as $n\rightarrow \infty$, tends to $e^{-\lambda}$.
Finally, as a special case, it is shown that, for fixed $p$ if $p=1/2$, the
probability that $Q_{p_n}^n$ is connected , as $n \rightarrow \infty$,
tend to $e^{-1}$; and if $p>1/2$, the probability that $Q_{p_n}^n$ is
connected, as $n\rightarrow \infty$, tends to $1$.

\section{Isolated vertices}
\noindent{\textbf{For $p<0.5$:}}\\
Assume $p$ is fixed and $p<0.5$. First, consider the following
definitions:
\begin{definition} $f_n(p_n)$:=\textbf{P}($Q_{p_n}^n$ is connected)\end{definition}

\begin{definition}$g_n(p_n):=\textbf{P}(Q_{p_n}^n\text{contains an isolated point})$
\end{definition}

\begin{definition}\label{1.d1}
$X_i(n):=\left\{
\begin{array}{ll}
         1 & \mbox{Vertex $i$ is isolated,  $i\in V(Q_{p_n}^n$)};\\
        0 & \mbox{Vertex $i$ is NOT isolated,  $i\in V(Q_{p_n}^n$)}.\end{array} \right.$
, and $X(n):=\displaystyle\sum_{i\in V(Q_{p_n}^n)}X_i(n)$.
\end{definition}

\noindent Now, calculate $E[X(n)]$ and $Var[X(n)]$ as follows:
\begin{align}
\mu: = E[X(n)] = \sum_{i\in V(Q_{p_n}^n)}E[X_i(n)] = \sum_{i\in
V(Q_{p_n}^n)}(1-p)^n=2^n(1-p)^n
\end{align}
\begin{align}
Var[X(n)] = \displaystyle\sum_{i\in V(
Q_{p_n}^n)}Var[X_i(n)]+\displaystyle\sum_{i\neq j; i,j\in V(Q_{p_n}^n)}
Cov[X_i(n),X_j(n)]
\end{align}
\noindent where, $Var[X_i(n)]$ and $Cov[X_i(n),X_j(n)]$ are equal
to:
\begin{align}
\displaystyle\sum_{i\in V(Q_{p_n}^n)}Var[X_i(n)]& =
2^n(1-p)^{n}-2^n(1-p)^{n}(1-p)^{n}=\mu-\mu(1-p)^{n}
\end{align}
\begin{align}
Cov[X_i(n),X_j(n)]& = E[X_i(n)X_j(n)]-E[X_i(n)]E[X_j(n)]\\
&= \left\{
\begin{array}{ll}
         0 & \mbox{i,j not adjacent};\\
        (1-p)^n(1-p)^{n-1}-(1-p)^n(1-p)^n=\frac{\mu^2}{2^{2n}}(\frac{p}{1-p})& \mbox{i,j adjacent}.\end{array} \right.
\end{align}

\noindent and finally:
\begin{align}
Var[X(n)]=\mu-\mu
(1-p)^n+\frac{\mu^2}{2^n}(\frac{np}{1-p})=\mu+\mu(1-p)^n(\frac{np}{1-p}-1)
\end{align}
Now, since we have $Var[X(n)]$, we can use Chebyshev's inequality to
estimate $g_n(p)$. Chebyshev's inequality states that:
\begin{align}
1-g_n(p)=\textbf{P}[X(n)=0]\leq \textbf{P}[|X(n)-\mu|\geq \mu]\leq
\frac{Var[X(n)]}{\mu^2}
\end{align}
By applying Chebyshev's inequality when $p<0.5$, one gets
$Var[X(n)]/\mu^2\rightarrow 0$, as $n\rightarrow\infty$. Therefore
$\lim_{n\rightarrow \infty} g_n(p)=1$. And finally, since
$f_n(p)\leq1-g_n(p)$, then for $p<0.5$ the probability that
$Q_{p_n}^n$ is connected for $p<0.5$, as $n\rightarrow \infty$, tends
to
$0$.\done\\\\

\noindent{\textbf{For $p>0.5$:}}\\
Assume $p$ is fixed and $p>0.5$. In order to calculate $g_n(p)$
when $p>0.5$, as $n\rightarrow \infty$, one can use the following
inequality:
\begin{align}
g_n(p)=\textbf{P}[X(n)>0]\leq E[X(n)]=\mu
\end{align}
Since $E[X(n)]\rightarrow 0$ as $n\rightarrow\infty$, then
$\lim_{n\rightarrow \infty}g_n(p)=0$. This means that the
probability that there are no isolated points in
$Q_{p_n}^n$ for $p>0.5$, as $n\rightarrow \infty$, tends to $1$.\done\\

\noindent{\textbf{For $p_n=1-(1/2)\lambda^{1/n}(1+o(1/n))$:}}\\
One needs the following theorem from \cite{bol} to find the distribution of $X(n)$ (distribution of
the number of isolated points).
\begin{theorem}\label{t.2.1}
Let $\lambda = \lambda(n)$ be a non-negative bounded function on
\textbf{N}. Suppose the non-negative integer valued random
variables $X(1),X(2),...$ are such that:
\begin{align}
\lim_{n\rightarrow \infty }\{E_r[X(n)]-\lambda^r\}=0, \text{ }
r=0,1,...
\end{align}
\noindent where $E_r[X]$ is the $r$th factorial moment of $X$,
i.e. $E_r[X]=E[(X)_r]$. Then
\begin{align}
X(n) \stackrel{d}{\longrightarrow}\textbf{P}_\lambda
\end{align}
\end{theorem}

\noindent Use the definition of $X(n)$ presented in definition
\ref{1.d1}. The goal is to calculate $E[X(n)]$.
\begin{align}
E_r[X(n)]=E[X(n)(X(n)-1)(X(n)-2)...(X(n)-r+1)]
\end{align}
Since $X(n):=\sum_{i\in V(Q_{p_n}^n)}X_i(n)$ and $X_i$'s are indicator
functions, therefore:
\begin{align}
X(n)(X(n)-1)(X(n)-2)...(X(n)-r+1)=\sum_{(i_1,i_2,...,i_r)}X_{i_1}X_{i_2}...X_{i_r}\label{1.31}
\end{align}
where the sum is over all ordered sets of distinct
vertices. Then:
\begin{align}
E_r[X(n)]&=E[X(n)(X(n)-1)(X(n)-2)...(X(n)-r+1)]\\
&=E[\sum_{(i_1,i_2,...,i_r)}X_{i_1}X_{i_2}...X_{i_r}]\\
&=\sum_{(i_1,i_2,...,i_r)}\textbf{P}[X_{i_1}=1,X_{i_2}=1,...,X_{i_r}=1]\label{1.32}
\end{align}
One knows that a set of $r$ vertices is incident with at most $rn$
edges. There are $(r)_r\binom{2^n}{r}$ ways to choose such $r$
vertices. Hence:
\begin{align}
 E_r[X(n)]\geq(r)_r\binom{2^n}{r}(1-p_n)^{rn}=(2^n)_r(1-p_n)^{rn}\label{1.26}
\end{align}
One the other hand, a set of $r$ vertices is incident with at
least $r(n-r)$ edges. There are at most
$(r-1)_{r-1}\binom{2^n}{r-1}(r-1)n$ ways to choose a set of $r$
vertices in $Q_{p_n}^n$ where at least two vertices are adjacent;
since if we choose $r-1$ vertices independently, then the last
vertex must be connected to one of the chosen vertices. In other
words, there are at most $(r-1)_{r-1}\binom{2^n}{r-1}(r-1)n$ ways
to choose $r$ vertices which some of them are adjacent to each
other. Hence:
\begin{align}
E_r[X(n)]&\leq
(2^n)_r(1-p_n)^{rn}+(r-1)_{r-1}\binom{2^n}{r-1}(r-1)n(1-p_n)^{r(n-r)}\\
&\leq
(2^n)_r(1-p_n)^{rn}+(2^n)_rrn(1-p_n)^{r(n-r)}\\
&\leq (2^n)_r(1-p_n)^{rn}+2^{n(r-1)}rn(1-p_n)^{r(n-r)}\label{1.29}
\end{align}
Finally from \ref{1.26} and \ref{1.29} one gets:
\begin{align}
(2^n)_r(1-p_n)^{rn}\leq E_r[X(n)]\leq (2^n)_r(1-p_n)^{rn}+2^{n(r-1)}rn(1-p_n)^{r(n-r)}\\
\end{align}
which gives:
\begin{align}
(2(1-p_n))^{rn}(1-\frac{r}{2^n})^r\leq
E_r[X(n)]\leq(2(1-p_n))^{rn}\{1+2^{-n}rn(1-p_n)^{-r^2}\}
\end{align}
Since $r$ is fixed and $\lim_{n\rightarrow
\infty}(2(1-p_n))^n=\lambda $, then:
\begin{align}
\lim_{n\rightarrow \infty}(E_r[X(n)])=\lambda^r \text{ for
r=0,1,2,...}
\end{align}
This shows that $ X(n)\stackrel{d}{\longrightarrow}\textbf{P}_\lambda$. \done\\\\

\noindent{\textbf{For $p=0.5$:}}\\
In the calculation of $p_n=1-1/2\lambda^{1/n}(1+o(1/n))$, if we fix $p=1/2$ and let $\lambda=1$, then we get that
 the distribution of $X(n)$ , as $n\rightarrow \infty$, tends to a Poisson distribution with mean $1$. Therefore, one can conclude:
\begin{align}
\lim_{n\rightarrow \infty}(1-g_n(p))=\lim_{n\rightarrow
\infty}(\textbf{P}(X(n)=0))=e^{-1}
\end{align}
\noindent This shows that for $p=1/2$ the probability that $Q_{p_n}^n$
has no isolated point, as $n\rightarrow \infty$, tends to
$e^{-1}$. \done

\section{Isoperimetric problem for the cube}
One needs an inequality which relates the size of a subgraph of
$Q^n$ to the size of its boundary. This inequality will be applied
to prove that for fixed values of $p$ if $p\geq 0.5$, then the
probability that subgraphs of $Q^n$ do not have a component of
size larger than $2$ and smaller than $2^n$, as $n\rightarrow \infty$, tends to $1$. The
proof presented here is based on the proof presented in
\cite{bol}.

\noindent\begin{definition}The edge boundary $b_G(H)$, where $H$
is an induced subgraph of G, is the number of edges which joins
vertices in $H$ to the vertices in $G\backslash
H$.\end{definition}
\noindent\begin{definition}$b_G(m):=\min\{b_G(H), H \text{ is an
induced subgraph of } G, |V(H)|=m\}$.\end{definition}

The main task in this section is to calculate $b_{Q^n}(m)$. The
answer, loosely, is if $m=2^k$ for some $k<n$ then one should take
a k-dimensional sub-cube of $Q^n$ as $b_{Q^n}(H)$. If $2^k\leq
m<2^{k+1}$, for some $k<n$, then one should choose one side of a
$(k+1)-$cube and $m-2^k$ more vertices properly chosen in the
other half. Since $Q^n$ is n-regular and $H$ is an induced
subgraph of $G$ with $|V(H)|=m$, then:
\begin{align}
b_{Q^n}(H)& = mn-2e(H),   \text{where e(H) is the total number of edges in H.}\\
b_{Q^n}(m)& = mn-2e_n(m), \text{where }e_n(m)=\max\{ e(H):H \text{
induced subgraph of } Q^n, |V(H)|=m \}.
\end{align}
\begin{definition}
$h(i)$ := sum of digits in the binary expansion of $i$ and
$f(l,m):=\displaystyle\sum_{l\leq i<m}h(i)$
\end{definition}

\noindent\begin{lemma}\label{1.1} If $1\leq k\leq l$, then
$f(l,l+k)\geq f(0,k)+k$
\end{lemma}
\proof \noindent Let look at the binary expansion of a few
numbers:
\[
\begin{array}{ccccc}
Column&3&2&1&0\\
Bin\backslash Dec&2^3&2^2&2^1&2^0\\
  0 &   &   &   & 0  \\
  1 &   &   &   & 1  \\
  2 &   &   & 1 & 0  \\
  3 &   &   & 1 & 1  \\
  4 &   & 1 & 0 & 0 \\
  5 &   & 1 & 0 & 1 \\
  6 &   & 1 & 1 & 0 \\
  7 &   & 1 & 1 & 1 \\
  8 & 1 & 0 & 0 & 0 \\
  9 & 1 & 0 & 0 & 1 \\
 10 & 1 & 0 & 1 & 0 \\
 11 & 1 & 0 & 1 & 1
\end{array}
\]

From this representation, one can observe that column $i$ starts
with a block of $2^i$ zeros. Therefore, sum of jth digits of $k$
consecutive numbers is minimal if the first block of $0$'s is as
long as possible. Hence, one can conclude:
\begin{align}\label{1.2}
f(l,l+k)\geq f(0,k)
\end{align}

For every $i$ define $r$ such that $ 0\leq i \leq 2^r-1$. The
binary expansions of $i$ and $2^r-1-i$ are symmetric. This means
that, if there is a $1/0$ in an specific location of the binary
expansion of $i$ then there is a $0/1$ in the same location of the
binary expansion of $2^r-1-i$. Therefore,
\begin{align}
h(i)+h(2^r-1-i)=r \text{ for } 0\leq i \leq 2^r-1
\end{align}
\noindent Consequently, since:
\begin{align}
\displaystyle\sum_{l\leq i<l+k}h(i)+\displaystyle\sum_{2^r-l-k\leq
i <2^r-l}h(i)=rk
\end{align}
\noindent then:
\begin{align}\label{1.5}
f(l,l+k)+f(2^r-l-k,2^r-l)=rk, \text{if } l+k\leq 2^r
\end{align}

\noindent Let us prove lemma \ref{1.1} with the assumption
$k\leq2^r \leq l$ by using inequalities \ref{1.2} and \ref{1.5}.
This assumption means that the length of the
sequence in the binary expansion of $2^r+k$ and $2^r$ are equal.\\

\noindent With the same logic that one gets \ref{1.2}, one gets:
\begin{align}
f(l,l+k)\geq f(2^r,2^r+k) \text{ when } 2^r\leq l
\end{align}
\noindent and then for $k \leq 2^r$ one can get:

\begin{align} \label{1.3}
f(2^r,2^r+k)=\displaystyle\sum_{2^r\leq i<2^r+k}h(i)
\end{align}

\noindent $\sum_{2^r\leq i<2^r+k}h(i)$ is the sum over numbers
with the same length in their binary expansion's sequence. When
one removes the last digit in their binary expansion, the remain
is $f(0,k)$. Therefore, $\sum_{2^r\leq i<2^r+k}h(i)$ is equal to
$k$ $1$'s plus $f(0,k)$. Hence:
\begin{align} \label{1.4}
f(2^r,2^r+k)=k+f(0,k) \text{ when } k\leq2^r \leq l
\end{align}
and finally:
\begin{align} \label{1.24}
f(l,l+k)\geq f(0,k)+k \text{ where }k\leq2^r\leq l
\end{align}

\noindent Now, one can prove lemma \ref{1.1} by induction on $K$,
without the assumption $k\leq2^r \leq l$. We want to prove that
for $1\leq K\leq l$, $f(l+K,l) \geq K + f(0,K)$. Fix $k$ such that
$1\leq k \leq l$ and $K<k$. For $K=1$ the inequality in lemma
\ref{1.1} is trivial. Assume that the inequality is true for $K<k$
and $K>2$, which means:
\begin{align} \label{1.7}
f(l,l+K)\geq K+f(0,K) \text{ when } 1\leq k \leq l \text{ , and }
K<k
\end{align}
Now, one should verify the inequality for $K=k$. Define $r\geq1$
by $2^{r-1}\leq k<2^r$. If $l\geq 2^r$, then $k\leq 2^r \leq l$ and the lemma is implied
by inequality \ref{1.24}. Hence, one may assume that $2^{r-1}< l <
2^r$. Now, one should apply inequality \ref{1.5} and \ref{1.7} in
order to get the final result:
\begin{align}
f(l+k)&=f(l,2^r)+f(2^r,l+k) \text{\emph{ (from definition of f and }} l\geq 2^r)\\
&=(2^r-l)r-f(0,2^r-l)+f(2^r,l+k) \text{ \emph{ (from \ref{1.5})}} \\
&\geq (2^r-l)r-f(0,2^r-l) +f(0,l+k-2^r)+l+k-2^r \text{\emph{ (from \ref{1.7})}} \\
&\geq (2^r-l)r-f(2^r-k,2^r-k+2^r-l)+2^r-l +f(0,l+k-2^r)+l+k-2^r \text{ \emph{(from \ref{1.7})}} \\
&\geq (2^r-l)r-f(2^r-k,2^r-k+2^r-l)+f(0,l+k-2^r)+k  \\
&\geq f(l+k-2^r,k)+f(0,l+k-2^r)+k \text{ \emph{(from \ref{1.5})}} \\
&\geq f(0,k)+k \text{ \emph{(from characteristics of f)}}
\end{align}
\done
\begin{theorem}\label{t.2}
For $2\leq m \leq 2^{n}$ we have $b_{Q^n}(m)=mn-2f(0,m)$. In
other words, $f(0,m)=e_n(m) \text{ where }e_n(m)=\max\{ e(H):H
\text{ induced subgraph of } Q^n, |V(H)|=m \}.$
\end{theorem}
\proof First, let us fix an $m$. As the first step one should
prove that $e_{n}(m)\geq f(0,m)$. Vertex $i$ is connected to
$h(i)$ vertices $j$ with $j<i$, since for each 1 in the binary
expansion of $i$ there is exactly one j $(j<i)$, which its binary
expansion differs in the position of that $1$. Therefore, one can
conclude that $W=\{0,1,2,...,m-1 \}$ contains $\sum_{0\leq
i<m}h(i)=f(0,m)$ edges. So, $e_{n}(m)\geq f(0,m)$.

As the second step, one should prove that $e_n(m)\leq f(0,m)$ by
induction on $n$. Fix $m$ and $n$ for $2\leq m \leq 2^{n}$. For
$n=1$ the inequality is trivially true. Assume that it is true for
$N<n$, which means:
\begin{align}
e_N(m)\leq f(0,m), \text{ where } N<n \text{ and the fixed m is: }
2\leq m \leq 2^{n}\label{1.21}
\end{align}
\noindent Now, one should check the inequality \ref{1.21} for
$N=n$. This means that we should find an $H$ induced subgraph of
$Q^n$, $|V(H)|=m$, which maximize $e_n(m)$. Let us split $Q^n$
into two (n-1)-dimensional cubes, the top face with $2^{n-1}$
vertices and the bottom face with $2^{n-1}$ vertices. This means,
there are $(n-1)2^{n-2}$ edges in each face, and $2^{n-1}$ edges
between two faces. Now, one can construct $H$. Choose $m_1$
vertices for $H$ from the top face, and $m_2$ vertices from the
bottom face, where $m_1+m_2=m$ and $m_1\leq m_2$. In other words,
$H$ is constructed from two induced subgraphs, one from the top
face, denoted by $H_1$, and the other from the bottom face,
denoted by $H_2$.

 Each face is a (n-1)-dimensional cube, so inequality \ref{1.21}
 holds for both $H_1$ and $H_2$. Also, each vertex of the top
face is connected to exactly one vertex of the bottom face. Hence,
the number of edges of $H$ is at most:
\begin{align}
e_{n}(m) \leq f(0,m_1) + f(0,m_2) + m_1 \text{(from \ref{1.21})}
\end{align}
\noindent where $m_1$, in the right hand side of the inequality,
is for the maximum number of edges between $H_1$ and $H_2$,
which one can choose here. Finally, by applying lemma \ref{1.1},
one gets:
\begin{align}
e_{n}(m) & \leq f(0,m_1) + f(0,m_2) + m_1\\
& \leq f(m_2,m_2+m_1)+f(0,m_2) \text{ (from lemma \ref{1.1})} \\
& \leq f(0,m) \text{ (from definition of f)}
\end{align}
\done\\
 Theorem \ref{t.2} shows that, if we want to choose an
induced subgraph of $Q^n$, with $m$ vertices, which has the
smallest edge boundary, then we should choose the induced subgraph
of $Q^n$ with the set of vertices $W=\{0,1,2,...,m-1\}$.\\
\begin{corollary}\label{1.6}
For all $k$ and $n$, $e_n(k)\leq \frac{k}{2}\lceil log_2k\rceil$,
which is equivalent to $b_{Q^n}(k)\geq k(n-\lceil log_2k\rceil)$.
\end{corollary}
\proof \noindent Let $r=\lceil log_2k\rceil$. Then
\begin{align}
2f(0,k)&\leq f(0,k)+f(0,k)\leq f(0,k)+f(2^r-k,2^r) \text{ (from \ref{1.2}) }\\
&= rk \text{ (from \ref{1.5})}
\end{align}
\noindent Therefore:
\begin{align}
e_{n}(k)=f(0,k)\leq r\frac{k}{2}=\frac{k}{2} \lceil log_2k
\rceil
\end{align}
\done

\section{Isolated components of size larger than 2 and smaller than $2^n$}
\begin{definition}
$C_s$ is the family of s-subsets (subsets with size s) of
$V=V(Q^n)$ whose induced graph is connected.
\end{definition}

\noindent\textbf{Remarks:} $h(n):=o(g(n))$ means $\frac{h(n)}{g(n)}\rightarrow 0$ as $n\rightarrow \infty$.\\
\noindent\textbf{Remarks:} The following inequality will be
applied a lot in the rest of this section:
\begin{align}
(\frac{n}{k})^k \leq \binom{n}{k} \leq \frac{n^k}{k!} \leq
(\frac{ne}{k})^k \label{rem.1}
\end{align}

\begin{theorem}
If $p_n\geq 1- \frac{1}{2}(\log n)^{\frac{1}{n}}$, the probability
that for some $S\in C_s$, $2\leq s \leq 2^{n-1}$, no edges
of $Q_{p_n}^n$ join $S$ to $V(Q^n)\setminus S$, as $n \rightarrow \infty$,
tends to 0.
\end{theorem}
\noindent \textbf{Note:} For $2^{n-1}<s<2^n$, if there exist a
component of size smaller than $2^n$ then there is at least one
component of size smaller than $2^{n-1}$ which contradicts with
the theorem.

\proof \noindent Consider $S\subset V=V(Q^n)$ and set
$b(S)=b_{Q^n}(H)$ for which $H$ is the induced subgraph of $Q^n$
with the set of vertices S. One can observe that:
\begin{align}
\textbf{P}(\text{No edges of }
Q_{p_n}^n \text{ join S to } V\setminus S)=
(1-p_n)^{b(S)}
\end{align}
\noindent In order to prove the theorem, it is sufficient to show
that:
\begin{align}
\sum_{s=2}^{2^{n-1}}\sum_{S\in C_s}(1-p_n)^{b(S)}=o(1)
\end{align}
From corollary \ref{1.6}, one knows that for $|S|=s$:
\begin{align}
b(S)\geq b(s)\geq s(n-\lceil \log_2s \rceil)
\end{align}
\noindent and therefore:
\begin{align}
\sum_{S\in C_s}(1-p_n)^{b(S)}\leq |C_s|(1-p_n)^{b(s)}
\end{align}
\noindent Now, one should partition s, $2\leq s \leq 2^{n-1}$, to
different intervals in order to find a bound for $|C_s|$ and
$(1-p_n)^{b(s)}$ for each
interval.\\\\
\noindent \textbf{First interval $2 \leq s \leq s_1, s_1=\lfloor \frac{2^{\frac{n}{2}}}{n^2} \rfloor$:}\\

First, one should find a bound for $|C_s|$. One has maximum $2^n$
choices to choose the first element for $C_s$. The selected element
is connected to maximum $n$ vertices, therefore there are $n$
choices to choose the second element. With the same logic there
are at most $(s-1)n$ choices to choose the last element for $C_s$.
Therefore, one can show:
\begin{align}
|C_s|\leq 2^n(n)(2n)...((s-1)n)) \leq (s-1)!(n)^{s-1}2^n
\end{align}
\noindent Hence:
\begin{align}
|C_s|(1-p_n)^{b(s)}\leq (s-1)!(n)^{s-1}2^n (1-p_n)^{s(n-\lceil
\log_2s \rceil)} \label{1.eq.2}
\end{align}
Since $p_n=1-\frac{1}{2}(\log n)^{\frac{1}{n}}$, so for large enough $n$:
\begin{align}
(1-p_n)^{s(n-\lceil log_2s \rceil)}&\leq(2)^{-ns}(\log n)^{s}(1-p_n)^{-s(\log_2s)} \text{(neglecting some small terms)}\\
& = (2)^{-ns}(\log n)^{s}2^{s\log_2s}(\log n)^{\frac{-s\log_2 s}{n}}\\
&(\text{ since for large enough n: } (\log n)^{\frac{-s\log_2s}{n}}\leq 1 )\\
&\leq (2)^{-ns}(\log n)^{s}s^s \label{1.eq.1}
\end{align}
\noindent From equations \ref{1.eq.2} and \ref{1.eq.1}, one can show
that:
\begin{align}
|C_s|(1-p_n)^{b(s)}\leq (s-1)!(n)^{s-1}2^n(2)^{-ns}(\log n)^{s}s^s
\label{1.eq.3}
\end{align}
\noindent Assume that the right hand side of inequality
\ref{1.eq.3} is equal to A. After multiplying both sides of
inequality \ref{1.eq.3} with $\frac{ns^{s+1}}{s!}$ and then taking
$\log_2$ from both sides, one gets:
\begin{align}
\log_2(|C_s|(1-p_n)^{b(s)}\frac{ns^{s+1}}{s!})\leq\log_2(A\frac{ns^{s+1}}{s!})
\end{align}
\noindent If $\log_2(A\frac{ns^{s+1}}{s!})\rightarrow-\infty$ as
$n\rightarrow\infty$ then $A\frac{ns^{s+1}}{s!}$ should tend to 0. This
means that $|C_s|(1-p_n)^{b(s)}\frac{ns^{s+1}}{s!}$ tends to 0, as
$n\rightarrow \infty$. Therefore:
\begin{align}
|C_s|(1-p_n)^{b(s)}\leq \frac{s!}{ns^{s+1}} \text{ for large
enough n}
\end{align}
which shows that:
\begin{align}
\sum_{s=2}^{s_1}\sum_{S\in C_s}(1-p_n)^{b(S)}=o(1)
\end{align}

\noindent Finally, it remains to prove $\log_2(A\frac{ns^{s+1}}{s!})\rightarrow-\infty$ as $n\rightarrow\infty$. One can verify this for $s\leq n$ and
$s>n$.\\\\

\noindent \textbf{Second interval $s_1 +1 \leq s \leq 2^{n-1}$ and
$S\in C_s^-, s_1= \lfloor \frac{2^{\frac{n}{2}}}{n^2} \rfloor$:}\\
Let us define $C_s^-$ and $C_s^+$ as follows:
\begin{definition}
\begin{align}
C_s^-:=\{S\in C_s |b(s) \geq s(n-\log_2s+\log_2n)\} \text{, and }
C_s^+:=C_s\backslash C_s^-
\end{align}
\end{definition}
One can bound $|C_s^-|$ for  $s_1 +1 \leq s \leq 2^{n-1}$  as
follows:
\begin{align}
|C_s^-|\leq |C_s| \leq \binom{2^n}{s} \leq \frac{2^{ns}}{s!} \leq
(\frac{e2^n}{s})^s
\end{align}
Hence:
\begin{align}
\displaystyle\sum_{s=s_1+1}^{2^{n-1}}\sum_{S\in C_s^-}(1-p_n)^{b(S)} &\leq \sum_{s=s_1+1}^{2^{n-1}}(\frac{e2^n}{s})^s(\frac{1}{2}(\log n)^{\frac{1}{n}})^{s(n-\log_2s+\log_2n)}\\
\displaystyle& \leq \sum_{s=s_1+1}^{2^{n-1}}(\frac{e2^n2^{-(n-\log_2 s +\log_2 n)}(\log n)^{\frac{(n-\log_2s+\log_2n)}{n}}}{s})^s\\
& \leq \sum_{s=s_1+1}^{2^{n-1}}(\frac{e2^n2^{-n}2^{\log_2s}2^{-\log_2n}\log n}{s})^s(\log n)^{\frac{n(-\log_2s+\log_2n)}{s}}\\
&(\text{ since for large enough n: } (\log n)^{\frac{n(-\log_2s+\log_2n)}{s}} \leq 1)\\
& \leq \sum_{s=s_1+1}^{2^{n-1}} (\frac{e\log n}{n})^s = o(1)
\end{align}
\done\\\\
\noindent \textbf{Third interval $s_1 \leq s \leq s_2$, $s_1=
\lfloor
\frac{2^{\frac{n}{2}}}{n^2} \rfloor, s_2=\lfloor \frac{2^n}{(\log n)^4}\rfloor$ and $S\in C_s^+$:}\\
For the 3rd and the 4th intervals one needs to know how to find a
bound for $|C_s^+|$. The following lemma, presented by B.Bollobas
\cite{bol}, helps us in this matter:
\begin{lemma}\label{1.l.1}
Let G be a graph of order $v$ and suppose that $\Delta(G)\leq
\Delta$, $2e(G)=vd$ and $\Delta +1 \leq u \leq v- \Delta -1$.
Then, there is a u-set of U of vertices with:
\begin{align}
|N(U)|=|U \cup \Gamma (U)| \geq v\frac{d}{\Delta} \{
1-exp(\frac{-u(\Delta+1)}{v}) \}
\end{align}
where, $\Delta(G):=$ Maximum degree in G, $d:=$ average degree in
G and $\Gamma(U)=\{x\in V(G): xy \in E(G) \text{ for some y}\in
U\}$
\end{lemma}
Let $H=Q_n[S]$ (the induced subgraph of $Q_n$ with the set of
vertices $S$). From the definition of $C_s^+$ one knows that the
average degree in $H$ is at least:
\begin{align}
\log_2s-\log_2n
\end{align}
The goal is to find $U\subset S$, where
$|U|:=u:=\lfloor\frac{2s}{n}\rfloor$, $\Delta = n$, $v=s$, $d\geq
\log_2s-\log_2n$ and then use lemma \ref{1.l.1} to calculate
$|N(U)|$. First, one should check the condition $\Delta +1 \leq u
\leq v- \Delta -1$ for defined variables in order to use lemma
\ref{1.l.1}. First, check if $n+1 \leq
\lfloor\frac{2s}{n}\rfloor$, as $n\rightarrow \infty$:
\begin{align}
&\frac{2s}{n}=\frac{2^{\frac{n}{2}+1}}{n^3} \text{for minimum s, and trivially }  n+1\leq \frac{2^{\frac{n}{2}+1}}{n^3} \text{ for large enough n} \\
\end{align}
\noindent and then check if $\lfloor\frac{2s}{n}\rfloor \leq
s-(n+1)$. One should check if  $ns-n(n+1)\geq 2s$, which means
one should check that whether:
\begin{align}\frac{2^{\frac{n}{2}}(n-2)}{n^3(n+1)}\geq1 \end{align}
\noindent which is clearly true for large enough $n$. Now, one can
apply lemma \ref{1.l.1} on the graphs generated by $S$ and get :
\begin{align}
\exists U\subset S : |N(U)|\geq
s\frac{\log_2s-\log_2n}{n}\{1-exp(-\frac{u(n+1)}{s})\}
\end{align}
\noindent where:
\begin{align}
\frac{\log_2s-\log_2n}{n}\geq \frac{(\log_2(\frac{2^{\frac{n}{2}}}{n^2})-\log_2n)}{n}=\frac{n-6\log_2n}{2n}\\
\text{and } \lim_{n\rightarrow\infty}\frac{n-6\log_2n}{2n}=
\frac{1}{2} \label{1.8}
\end{align}
\noindent on the other hand:
\begin{align}
\lim_{n\rightarrow\infty}(1-exp(-\frac{n+1}{s}(\frac{2s}{n}+1)))=1-e^{-2}
\label{1.9}
\end{align}
Therefore, from \ref{1.8} and \ref{1.9} one gets:
\begin{align}
|N(U)| \geq \frac{1}{2}(1-e^{-2})s \geq \frac{s}{3} \text{ as } n\rightarrow \infty\label{1.10}
\end{align}

\noindent Now that we have $|N(U)|$, we can estimate a bound for
$|C_s^+|$ here. We know from \ref{1.10} that for each $S \in
C_s^+$ there exist a $U \subseteq S$,
$|U|:=u:=\lfloor\frac{2s}{n}\rfloor$, such that
$|N(U)|\geq s/3$. Therefore, one can choose $S\in C_s^+$ as follows:\\
\noindent 1. Select $u$ vertices of $Q^n$; there are $\binom{2^n}{u}$ choices for this $u$.\\
\noindent 2. Select $\lfloor\frac{s}{3}\rfloor-u$ neighbors of the
selected vertices in part $1$; there are maximum $(2^{n})^u$
choices, since there are at most
$\binom{n}{0}+\binom{n}{1}+\binom{n}{2}+...\binom{n}{n}=2^{n}$
ways to find neighbors of a vertex in $U$.\\
\noindent 3. Select $\lfloor\frac{2s}{3}\rfloor$ other vertices; there are at most $\binom{2^n}{\lfloor\frac{2s}{3}\rfloor}$ choices.\\
\noindent Hence:
\begin{align}
|C_s^+|\leq
\binom{2^n}{u}(2^{n})^u\binom{2^n}{\lfloor\frac{2s}{3}\rfloor}
\end{align}
\noindent \noindent and:
\begin{align}
\sum_{S\in C_s^+}(1-p_n)^{b(S)}\leq
\binom{2^n}{u}(2^{n})^u\binom{2^n}{\lfloor\frac{2s}{3}\rfloor}(1-p_n)^{b(s)}
\label{1.11}
\end{align}
\noindent where:
\begin{align}
(1-p_n)^{b(s)} \leq 2^{-sn}s^s(\log n)^{s} \label{1.12}
\end{align}
\noindent consequently from \ref{1.11}, \ref{1.12} and
\ref{rem.1}:
\begin{align}
\sum_{S\in C_s^+}(1-p_n)^{b(S)}\leq
(\frac{e2^n}{u})^u2^{un}(\frac{e2^n}{\lfloor\frac{2s}{3}\rfloor})^{\lfloor\frac{2s}{3}\rfloor}2^{-sn}s^s(\log
n)^{s}\label{1.14}
\end{align}
\noindent Write $s=2^{\beta n}$, ($\beta = \frac{log_2s}{n}$), so
that:
\begin{align}
2^{\beta n} \leq \frac{2^n}{(\log n)^4}  \Rightarrow \beta \leq
1-\frac{4\log_2\log n}{n} \label{1.13}
\end{align}
\noindent Now, find a bound for the inequality \ref{1.14}. First
calculate the first part of the inequality:
\begin{align}
(\frac{e2^n}{u})^u 2^{un}
(\frac{e2^n}{\lfloor\frac{2s}{3}\rfloor})^{\lfloor\frac{2s}{3}\rfloor}&
\leq (\frac{e2^n}{\frac{2s}{n}})^{\frac{2s}{n}}2^{2s}
(\frac{e2^n}{\frac{2s}{3}})^{\frac{2s}{3}}
=(2^22^2(\frac{3}{2}e)^{\frac{2}{3}})^s\frac{2^{\frac{2s}{3n}}}{s^{\frac{2s}{3}}}(\frac{2^n}{s})^{\frac{2s}{3}}\\
&(\text{ since for large enough n and $s_1 \leq s \leq s_2$: } (\frac{e2^n}{\frac{2s}{3}})^{\frac{2s}{3}}\leq 1 )\\
&=(2^22^2(\frac{3}{2}e)^{\frac{2}{3}})^s\frac{2^{\frac{2s}{3n}}}{s^{\frac{2s}{3}}}
=c^s2^{\frac{2}{3}sn(1-\frac{\log_2s}{n})}=c^s2^{\frac{2}{3}sn(1-\beta)}\label{1.15}
\end{align}
\noindent where c is a positive constant. Now, by substituting
\ref{1.15} in \ref{1.14} one gets:
\begin{align}
\sum_{S\in C_s^+}(1-p_n)^{b(S)} & \leq 2^{-sn}s^s(\log n)^{s}c^s2^{\frac{2}{3}sn(1-\beta)} \\
&=c^s(\log n)^{s}2^{-\frac{sn(1-\beta)}{3}}\\
&\leq c^s(\log n)^{s}2^{-\frac{4s\log_2\log n}{3n}} \text{ , (from \ref{1.13})} \\
&=c^s(\log n)^{s}2^{\log_2(\log n)^{\frac{-4s}{3}}}\\
&\leq c^s(\log n)^{s}(\log n)^{\frac{-4s}{3}}\\
& =c^s(\log n)^{\frac{-s}{3}}\label{1.25}
\end{align}
\noindent and finally from \ref{1.25}:
\begin{align}
\sum_{s=s_1}^{s_2}\sum_{S\in C_s^+}(1-p_n)^{b(S)}\leq
\sum_{s=s_1}^{s_2} c^s(\log n)^{\frac{-s}{3}} = o(1)
\end{align}
\done\\\\
\noindent \textbf{Fourth interval $s_2+1 \leq s \leq 2^{n-1}$ and $,s_2=\lfloor \frac{2^n}{(\log n)^4}\rfloor, S\in C_s^+$:}\\

In $H=Q^n[S]$ (the induced subgraph of $Q^n$ with the set of
vertices S), the average degree is at least:
\begin{align}
\log_2s-\log_2n>n-2\log_2n \label{1.33}
\end{align}
\noindent since:
\begin{align}
s\geq \lceil\frac{2^n}{(\log n)^4}\rceil \Rightarrow
\log_2(\frac{2^n}{(\log n)^4})<\log_2 s\\
\Rightarrow \log_2s - \log_2 n \geq \log_2(\frac{2^n}{(\log
n)^4})-\log_2n \geq n -\log_2(\log n)^4-\log_2n\\
(\text{for large enough n one can get, } n > (\log n)^4)\\
\geq n-2\log_2n
\end{align}
\noindent First, look for a subgraph of H with large average
degree. Let T be the set of vertices of H with degree at least
$n-(\log_2n)^2$, and set $t=|T|$. From \ref{1.33} one can conclude
that the sum of degrees in $H$ is at least $s(n-2\log_2n)$.We also
know that:
\begin{align}
\text{Sum of degrees in }S & \leq s(n-2\log_2n)\\
&\leq t\times(\text{Maximum degree of vertices in set $T$ of graph $H$ })\\
&+ (s-t)\times(\text{Maximum degree of vertices in set $S\setminus
T$ of
graph $H$ })\\
& \leq tn +(s-t)(n-(\log_2n)^2)\\
\Rightarrow t\geq s(1-\frac{2}{\log_2n})\label{1.34}
\end{align}

\noindent Define $H_1=Q^n[T]=H[T]$ as the induced subgraph spanned
by $T$. We want to calculate $|N_{H_1}(U)|$ for some $U$ in $H_1$,
hence we should estimate the size of $H_1$ and after that
calculate the average degree in T. Let us first calculate
$e(H_1)$, the total number of edges in $H_1$.
\begin{align}
e(H_1)\geq e(H)-(s-t)n
\geq\frac{s}{2}(n-2\log_2n)-\frac{2s}{\log_2n}n \text{ (from
\ref{1.33} and \ref{1.34})}
\end{align}
\noindent One knows that the average degree in $H_1$ is at least
$\frac{2e(H_1)}{s}$, and:
\begin{align}
\frac{2e(H_1)}{s}\geq n-2\log_2n-\frac{4n}{\log_2n}\geq
n-\frac{5}{log_2n}\\
(\text{since: } \log_2n^2<\frac{n}{\log_2n} \text{ for large
enough n})
\end{align}
\noindent Set $u=\lfloor\frac{2^n}{n^{\frac{1}{2}}}\rfloor$. One
should check the conditions of lemma \ref{1.1} here. Let
$v=t,\Delta=n, d\geq n-\frac{5}{\log_2n}$. So, one should check if
$n+1\leq\frac{2^n}{n^{\frac{1}{2}}}\leq t-(n+1)$ for large enough
$n$. Clearly, $n+1\leq\frac{2^n}{n^{\frac{1}{2}}}$, as
$n\rightarrow \infty$. It remains to prove
$\frac{2^n}{n^{\frac{1}{2}}}\leq t-(n+1)$, for large enough $n$.
For minimum $s$ from \ref{1.34} we can get:
\begin{align}
t\geq& \frac{2^n}{(\log n)^4}(1-\frac{2}{\log_2n}) \text{ (from \ref{1.34})}\\
&\geq \frac{2^n}{n^\frac{1}{2}} +n+1\text{ (for large enough n) }
\end{align}
Now, one can use lemma \ref{1.1} and estimate $|N_{H_1}(U)|$.
\begin{align}
|N_{H_1}(U)|&\geq
    \frac{t}{n}(n-\frac{5n}{\log_2n})\{1-\exp(-\frac{n+1}{t}\frac{2^n}{n^{\frac{1}{2}}})\}\\
&\geq\frac{t}{2}(1-\frac{5}{\log_2n})\{1-\exp(-\frac{n+1}{t}\frac{2^n}{n^\frac{1}{2}})\}
\label{1.18}
\end{align}
\noindent After that, let us estimate a bound for
$\exp(-\frac{n+1}{t}\frac{2^n}{n^{\frac{1}{2}}})$. One knows that
$t\geq s(1-\frac{2}{\log_2n})$. Since $\max(t)=s$ and
$\max(s)=2^{n-1}$, then:
\begin{align}
\frac{2^n(n+1)}{n^{\frac{1}{2}}t}\geq
\frac{2^n(n+1)}{n^{\frac{1}{2}}2^{n-1}} =
\frac{2(n+1)}{n^{\frac{1}{2}}}\geq n^\frac{1}{4} \text{( for large enough n)}\\
\Rightarrow \{1-\exp(-\frac{n+1}{t}\frac{2^n}{n^\frac{1}{2}})\}
\geq \exp(-n^\frac{1}{4})\text{( for large enough n)} \label{1.19}
\end{align}
\noindent By using the bound from \ref{1.19} in \ref{1.18}, one
gets:
\begin{align}
|N_{H}(U)|\geq|N_{H_1}(U)|& \geq t(1-\frac{5}{\log_2n})\{
1-\exp(-n^{\frac{1}{4}})\}\\
&=t\{1+\exp(-n^{\frac{1}{4}})\frac{5}{\log_2n}-\exp(-n^{\frac{1}{4}})-\frac{5}{\log_2n}\}\\
&(\lim_{n\rightarrow\infty}\exp(-n^{\frac{1}{4}})\frac{5}{\log_2n}=0 \text{ and } \exp(-n^{\frac{1}{4}})<
\frac{1}{\log_2n} \text{( for large enough n)})\\
&\geq
t\{1-\frac{1}{\log_2n}-\frac{6}{\log_2n}\}=t(1-\frac{6}{\log_2n})\\
&\geq s(1-\frac{2}{\log_2n})(1-\frac{5}{\log_2n})=s(1+\frac{2}{\log_2n}\frac{6}{\log_2n}-\frac{8}{\log_2n}) \text{(from \ref{1.34})}\\
&\geq s(1-\frac{8}{\log_2n})\label{1.20}
\end{align}

\noindent Now that we have $|N_{H}(U)|$, we can estimate a bound
for $|C_s^+|$ here. We know from \ref{1.20} that for each $S \in
C_s^+$ there exist a $U \subseteq S$,
$|U|:=u:=\lfloor\frac{2^n}{n^\frac{1}{2}}\rfloor$, such that
$|N_{H}(U)|\geq s(1-\frac{8}{\log_2n})$. Therefore, one can choose $S\in C_s^+$ as follows:\\
\noindent 1. Select $u$ vertices of $Q^n$; there are $\binom{2^n}{u}$ choices for this $u$.\\
\noindent 2. Select $\lfloor s(1-\frac{8}{\log_2n})\rfloor-u$
neighbors of the selected vertices in part 1. At most
$(\log_2n)^2$ of the $n$ neighbors of a vertex in $U$ do not
belong to $N_{H}(U)$. Hence there are at most
 $\sum_{(k_j)}(\prod_{i=1}^{u}\binom{n}{j})$ ways to find
neighbors of $u$ vertices in $U$, where the sum is over all
$(k_1,k_2,...,k_u), k_i \leq (\log_2n)^2$. We know that:
\begin{align}
\sum_{(k_i)}\prod_{i=1}^{u}\binom{n}{k_i}  &\leq \sum_{(k_i)}\prod_{i=1}^{u}(\frac{n^i}{i!})\\
&\leq \sum_{(k_i)}\prod_{i=1}^{u}(\frac{n^{(\log_2n)^2}}{(\log_2n)^2!})\\
&= \sum_{(k_i)}\frac{n^{u(\log_2n)^2}}{((\log_2n)^2!)^u}\\
&= {((\log_2n)^2)}^u\frac{n^{u(\log_2n)^2}}{((\log_2n)^2!)^u}\\
& \leq n^{u(\log_2n)^2}\\
\end{align}
\noindent 3. Select $\lfloor\frac{8s}{\log_2n}\rfloor$ other vertices; there are at most $\binom{2^n}{\lfloor\frac{8s}{\log_2n}\rfloor}$ choices.\\
\noindent Hence:
\begin{align}
\sum_{S\in C_s^+}(1-p_n)^{b(S)}\leq \binom{2^n}{u}n^{
u(\log_2n)^2}\binom{2^n}{\lfloor\frac{8s}{\log_2n}\rfloor}2^{-s(n-\log_2s)}(\log
n)^{s(1-\frac{\log_2s}{n})}\\
\end{align}
\noindent where:
\begin{align}
\binom{2^n}{u}n^{ u(\log_2n)^2}\binom{2^n}{\lfloor
\frac{8s}{\log_2n}\rfloor} \leq  2^{o(s)}
\end{align}
\noindent Therefore:
\begin{align}
\sum_{S\in C_s^+}(1-p_n)^{b(S)}\leq 2^{\varepsilon(s)}
\end{align}
\noindent where:
\begin{align}
\varepsilon(s)=o(s)-s\{n-\log_2s-\log_2\log n
+\frac{\log_2s}{n}\log_2\log n\}
\end{align}
\noindent Since $s\leq 2^{n-1}$, hence one can get:
\begin{align}
\varepsilon(s)&\leq o(s)-s\{n-(n-1)-\log_2\log n
+\frac{n-1}{n}\log_2\log n\}\\
&=o(s)-s\{1-\frac{1}{n}\log_2\log n\}\leq\frac{-s}{2}
\end{align}
\noindent Therefore, for large enough n, one can get:
\begin{align}
\sum_{s=s_2+1}^{2^{n-1}}\sum_{S\in C_s^+}(1-p_n)^{b(S)}\leq
\sum_{s=s_2+1}^{2^{n-1}}2^{-\frac{s}{2}}=o(1)
\end{align}
\done\\
\section{Comparison between a typical random graph and $Q^n$}
Consider $Per(G(V,E),p_n)$ as the subgraph of $G(V,E)$ after
percolation with the parameter $p_n$. The goal in this section is
to compare $\textbf{P}(Q^n_{p_n} \text{ is connected})$ with
$\textbf{P}$($Per(G\in G(N=2^n,M=n2^{n-1}),p_n$) is connected) as
$n\rightarrow \infty$.

One knows, the probability that $G_p \in G(n,p)$ is connected, for
$p=c \log n/n$ as $n\rightarrow \infty$, tends to \cite{bol} :
\begin{align}
\lim_{n\rightarrow \infty} \textbf{P}(G_p\in G(n,p) \text{ is
connected})= \left\{
\begin{array}{ll}
         1 & \mbox{if $c > 1$};\\
         1 -e^{-1}& \mbox{if $c=1$};\\
        0 & \mbox{if $c < 1$}.\end{array} \right.
\end{align}

\begin{theorem} If $Q$ is a convex property and $pq\binom{n}{2}\rightarrow
\infty$, then almost every graph in $G(n,p)$ has $Q$ iff for every
fixed $x$ a.e. graph in $G(n,M)$ has $Q$, when $M=\lfloor
p\binom{n}{2}+x(pq\binom{n}{2}^{0.5})\rfloor$ \cite{bol}.
\label{t6.1}
\end{theorem}
\begin{definition} $Q$ is a convex property if $F\subset G \subset
H$ and $F\in Q$ and $H \in Q$ then $G\in Q$.
\end{definition}

 We want to calculate $\textbf{P}(Per(G_m\in G(n,M=m),p_n) \text{ is
 connected})$ as $n\rightarrow \infty$. From theorem \ref{t6.1}, one can conclude that $G(n,p)$ and
$G(n,M=p\binom{n}{2})$ have almost the same behavior for
connectivity, as $n\rightarrow \infty$. Hence:
\begin{align}
\textbf{P}(Per(G_M\in G(n,M=p\binom{n}{2}),p_n) \text{ is
 connected})\\
 \approx \textbf{P}(G_p \in G(n,pp_n)) \text{ is
 connected}), \text{ as }n\rightarrow \infty
\end{align}
We can calculate $\textbf{P}(G_p \in G(n,pp_n)) \text{ is
 connected})$ from theorem \ref{t6.1} for $pp_n=c \log n/n$, as $n\rightarrow
 \infty$. Therefore, we should consider $p_n= c(n-1)\log n /(2m)$.
 Finally, one can calculate the probability that $G_M\in G(n,M=m)$ is connected after
percolation with the parameter $p_n= c(n-1)\log n /(2m)$ as:
\begin{align}
\lim_{n\rightarrow \infty} \textbf{P}(Per(G_M\in G(n,M),p_n)
\text{ is connected})= \left\{
\begin{array}{ll}
         1 & \mbox{for $p_n$'s such that $c > 1$};\\
         1 -e^{-1}& \mbox{for $p_n$'s such that $c=1$};\\
        0 & \mbox{for $p_n$'s such that $c < 1$}.\end{array} \right.\label{1.35}
\end{align}
Now, let us calculate $\textbf{P}(Per(G_M\in G(2^n,n2^{n-1}),p_n)
\text{ is connected)}$. From \eqref{1.35}, for
$p_n=c(1-\frac{1}{2^n})$, one gets:
\begin{align}
\lim_{n\rightarrow \infty} \textbf{P}(Per(G_M\in
G(2^n,n2^{n-1}),p_n) \text{ is connected})= \left\{
\begin{array}{ll}
         1 & \mbox{for $p_n$'s such that $c > \log2$};\\
         1 -e^{-1}& \mbox{for $p_n$'s such that $c=\log2$};\\
        0 & \mbox{for $p_n$'s such that $c < \log2$}.\end{array} \right.\label{1.36}
\end{align}
Finally, it can be concluded that, $\lim_{n\rightarrow \infty}
\textbf{P}(Per(G_M\in G(2^n,n2^{n-1}),p') \text{ is connected})=
1$ for $0.3\approx\log2<p'<0.5$, while $\lim_{n\rightarrow \infty}
\textbf{P}(Q^n_{p'} \text{ is connected})=0$ for $\log2<p'<0.5$.
Hence, a typical graph with $2^n$ vertices and $n2^{n-1}$ edges is
"more connected" than $Q^n$.

\chapter{Connected random subgraphs of the 3-cube}
\footnote{The proof presented for the 3-cube, in this chapter, is independent of the proof presented in the previous
 chapter; therefore, some parts of the proofs overlap each other.}
The 3-cube graph, $^3Q^n$, is a simple graph with $3^n$ vertices.
If one labels each vertex of $^3Q^n$ from $0$ to $3^n-1$, then two
vertices are adjacent if their ternary representation differs only
in one digit. The number of edges in $^3Q_n$ is $n3^n$, since each
vertex is connected to $2n$ vertices. A random subgraph of $^3Q^n$
contains all vertices of $^3Q^n$, and each edge independently with
probability $p_n$. $^3Q_{p_n}^n$ stands for a random subgraph of
$^3Q^n$.

The main goal in this chapter is to explore a critical value
$p_c$, which for fixed values of $p$ if $p<p_c$ then the
probability that $^3Q_{p_n}^n$ is connected, as $n \rightarrow
\infty$, tends to $0$; but if $p>p_c$ then the probability that
$^3Q_{p_n}^n$ is connected, as $n \rightarrow \infty$, tends to $1$.
We prove that this critical value is $(\sqrt3-1)/(\sqrt3)$.

In the first section of this chapter, first the probability that
$^3Q_{p_n}^n$ has no isolated point as $n\rightarrow \infty$, for
fixed $p$, is investigated. It is proved that for
$p<(\sqrt3-1)/(\sqrt3)$  the probability that $^3Q_{p_n}^n$ has no
isolated point, as $n \rightarrow \infty$, tends to
 $0$. Therefore, for $p<(\sqrt3-1)/(\sqrt3)$ the probability that $^3Q_{p_n}^n$ is connected, as $n \rightarrow \infty$, tends to
 $0$. Then it is proved that, for $p>(\sqrt3-1)/(\sqrt3)$ the probability that $^3Q_{p_n}^n$ has no isolated
point, as $n\rightarrow \infty$,  tends to $1$. In the next step,
the probability that $^3Q_{p_n}^n$ has no isolated point, as
$n\rightarrow \infty$, when $p$ depends on $n$ and it is close
$(\sqrt3-1)/(\sqrt3)$, is explored. It is proved that for
$\lambda(n)=\lambda>0$ and
$p_n=1-(1/\sqrt3)\lambda^{1/2n}(1+o(1/n)$, the probability that
$^3Q_{p_n}^n$ has no isolated point, as $n\rightarrow \infty$, tends
to $e^{-\lambda}$. Finally, it is proved that, for fixed
$p=(\sqrt3-1)/(\sqrt3)$ the probability that $^3Q_{p_n}^n$ has no
isolated point, as $n\rightarrow \infty$, tends to $e^{-1}$.

In the second section, one sheds light on the Isoperimetric
problem, which is the problem of finding an inequality which
relates the size of a subgraph to the size of its boundary. One
needs such an inequality to explore the probability that $^3Q_{p_n}^n$
has a component which is not the whole graph.

In the last section, the Isoperimetric inequality is applied to
prove that when $p$ depends on $n$ and $p_n\geq1-(1/\sqrt3)(\log
n)^{1/n}$, then the probability that there are no components with
size larger than $1$ and smaller than $3^n$ in $^3Q_{p_n}^n$ , as
$n\rightarrow \infty$, tends to $1$. Therefore, for
$p_n=1-(1/\sqrt3)\lambda^{1/2n}(1+o(1/n))$, the probability that
$^3Q_{p_n}^n$ is connected, as $n\rightarrow \infty$, tends to
$e^{-\lambda}$. Finally, as a special case, it is shown that, for fixed $p$ if
$p=(\sqrt3-1)/(\sqrt3)$, the probability that $^3Q_{p_n}^n$ is
connected , as $n \rightarrow \infty$, tend to $e^{-1}$; and if
$p>(\sqrt3-1)/(\sqrt3)$, the probability that $^3Q_{p_n}^n$ is
connected, as $n\rightarrow \infty$, tends to $1$.

\section{Isolated vertices}
\noindent{\textbf{For $p_n<\frac{\sqrt3-1}{\sqrt3}$:}}\\\\
 First, one should consider the following definitions:
\begin{definition} $f_n(p_n)$:=\textbf{P}($^3Q_{p_n}^n$ is connected)\end{definition}

\begin{definition}$g_n(p_n):=\textbf{P}(^3Q_{p_n}^n\text{contains an isolated point})$
\end{definition}

\begin{definition}\label{2.d1}
$X_i(n):=\left\{
\begin{array}{ll}
         1 & \mbox{Vertex $i$ is isolated,  $i\in V(^3Q_{p_n}^n$)};\\
        0 & \mbox{Vertex $i$ is NOT isolated,  $i\in V(^3Q_{p_n}^n$)}.\end{array} \right.$
, and $X(n):=\displaystyle\sum_{i\in V(^3Q_{p_n}^n)}X_i(n)$.
\end{definition}

\noindent Now, calculate $E[X(n)]$ and $Var[X(n)]$ as follows:
\begin{align}
\mu := E[X(n)] = \sum_{i\in V(^3Q_{p_n}^n)}E[X_i(n)] = \sum_{i\in
V(^3Q_{p_n}^n)}(1-p)^{2n}=3^n(1-p)^{2n}
\end{align}
\begin{align}
Var[X(n)] = \displaystyle\sum_{i\in V(
^3Q_{p_n}^n)}Var[X_i(n)]+\displaystyle\sum_{i\neq j; i,j\in
V(^3Q_{p_n}^n)} Cov[X_i(n),X_j(n)]
\end{align}
\noindent where, $Var[X_i(n)]$ and $Cov[X_i(n),X_j(n)]$ are equal
to:
\begin{align}
\displaystyle\sum_{i\in V(^3Q_{p_n}^n)}Var[X_i(n)]& =
3^n(1-p)^{2n}-3^n(1-p)^{2n}(1-p)^{2n}=\mu-\mu(1-p)^{2n}
\end{align}
\begin{align}
Cov[X_i(n),X_j(n)]& = E[X_i(n)X_j(n)]-E[X_i(n)]E[X_j(n)]\\
&= \left\{
\begin{array}{ll}
         0 & \mbox{i,j not adjacent};\\
        (1-p)^{2n}(1-p)^{2n-1}-(1-p)^{2n}(1-p)^{2n}=\frac{\mu^2}{3^{2n}}(\frac{p}{1-p})& \mbox{i,j adjacent}.\end{array} \right.
\end{align}

\noindent and finally:
\begin{align}
Var[X(n)]=\mu-\mu
(1-p)^{2n}+\frac{\mu^2}{3^n}(\frac{2np}{1-p})=\mu+\mu(1-p)^{2n}(\frac{2np}{1-p}-1)
\end{align}
Now, since we have $Var[X(n)]$, we can use Chebyshev's inequality to
estimate $g_n(p)$. Chebyshev's inequality states that:
\begin{align}
1-g_n(p)=\textbf{P}[X(n)=0]\leq \textbf{P}[|X(n)-\mu|\geq \mu]\leq
\frac{Var[X(n)]}{\mu^2}
\end{align}
By applying Chebyshev's inequality when $p<(\sqrt3-1)/(\sqrt3)$,
one gets $Var[X(n)]/\mu^2\rightarrow 0$, as $n\rightarrow\infty$.
Therefore $\lim_{n\rightarrow \infty} g_n(p)=1$. And finally,
since $f_n(p)\leq1-g_n(p)$, then for $p<(\sqrt3-1)/(\sqrt3)$ the
probability that $^3Q_{p_n}^n$ is connected, as $n\rightarrow \infty$, tends to $0$.\done\\

\noindent{\textbf{For $p_n>\frac{\sqrt3-1}{\sqrt3}$:}}\\
 Assume $p$ is fixed and $p>(\sqrt3-1)/(\sqrt3)$. In order to calculate $g_n(p)$ when
$p>(\sqrt3-1)/(\sqrt3)$, as $n\rightarrow \infty$, one can use the
following inequality:
\begin{align}
g_n(p)=\textbf{P}[X(n)>0]\leq E[X(n)]=\mu
\end{align}
Since $E[X(n)]\rightarrow 0$ as $n\rightarrow\infty$, then
$\lim_{n\rightarrow \infty}g_n(p)=0$. This means that the
probability that there are no isolated points in
$^3Q_{p_n}^n$ for $p>(\sqrt3-1)/(\sqrt3)$, as $n\rightarrow \infty$, tends to
$1$.\done\\\\
\noindent{\textbf{For $p_n=1-(1/\sqrt3)\lambda^{1/2n}(1+o(1/n))$:}}\\
One needs the following theorem from \cite{bol} to find the
distribution of $X(n)$ (distribution of the number of isolated
points).
\begin{theorem}\label{t.2.1}
Let $\lambda = \lambda(n)$ be a non-negative bounded function on
\textbf{N}. Suppose the non-negative integer valued random
variables $X(1),X(2),...$ are such that:
\begin{align}
\lim_{n\rightarrow \infty }\{E_r[X(n)]-\lambda^r\}=0, \text{ }
r=0,1,...
\end{align}
\noindent where $E_r[X]$ is the $r$th factorial moment of $X$,
i.e. $E_r[X]=E[(X)_r]$. Then
\begin{align}
X(n) \stackrel{d}{\longrightarrow}\textbf{P}_\lambda
\end{align}
\end{theorem}

\noindent Use the definition of $X(n)$ presented in definition
\ref{2.d1}. The goal is to calculate $E[X(n)]$.
\begin{align}
E_r[X(n)]=E[X(n)(X(n)-1)(X(n)-2)...(X(n)-r+1)]
\end{align}

Since $X(n):=\sum_{i\in V(Q_p^n)}X_i(n)$ and $X_i$'s are indicator
functions, therefore:
\begin{align}
X(n)(X(n)-1)(X(n)-2)...(X(n)-r+1)=\sum_{(i_1,i_2,...,i_r)}X_{i_1}X_{i_2}...X_{i_r}\label{2.31n}
\end{align}
where the sum is over all ordered sets of distinct
vertices. Then:
\begin{align}
E_r[X(n)]&=E[X(n)(X(n)-1)(X(n)-2)...(X(n)-r+1)]\\
&=E[\sum_{(i_1,i_2,...,i_r)}X_{i_1}X_{i_2}...X_{i_r}]\\
&=\sum_{(i_1,i_2,...,i_r)}\textbf{P}[X_{i_1}=1,X_{i_2}=1,...,X_{i_r}=1]\label{2.32n}
\end{align}
One knows that a set of $r$ vertices is incident with at most
$2rn$ edges. There are $(r)_r\binom{3^n}{r}$ ways to choose such
$r$ vertices. Hence:
\begin{align}
 E_r[X(n)]\geq(r)_r\binom{3^n}{r}(1-p_n)^{2rn}=(3^n)_r(1-p_n)^{2rn}\label{2.26n}
\end{align}
One the other hand, a set of $r$ vertices is incident with at
least $r(2n-r)$ edges. There are at most
$(r-1)_{r-1}\binom{3^n}{r-1}(r-1)2n$ ways to choose a set of $r$
vertices in $^3Q_{p_n}^n$ where at least two vertices are
adjacent; since if we choose $r-1$ vertices independently, then
the last vertex must be connected to one of the chosen vertices.
In other words, there are at most
$(r-1)_{r-1}\binom{3^n}{r-1}(r-1)2n$ ways to choose $r$ vertices
which some of them are adjacent to each other. Hence:

\begin{align}
E_r[X(n)]&\leq
(3^n)_r(1-p_n)^{2rn}+(r-1)_{r-1}\binom{3^n}{r-1}2(r-1)n(1-p_n)^{r(2n-r)}\\
&\leq
(3^n)_r(1-p_n)^{2rn}+(3^n)_r2rn(1-p_n)^{r(2n-r)}\\
&\leq
(3^n)_r(1-p_n)^{2rn}+3^{n(r-1)}2rn(1-p_n)^{r(2n-r)}\label{2.29n}
\end{align}
Finally from \ref{2.26n} and \ref{2.29n} one gets:
\begin{align}
(3^n)_r(1-p_n)^{2rn}\leq E_r[X(n)]\leq (3^n)_r(1-p_n)^{2rn}+3^{n(r-1)}2rn(1-p_n)^{r(2n-r)}\\
\end{align}
which gives:
\begin{align}
(3(1-p_n)^2)^{rn}(1-\frac{r}{3^n})^r\leq
E_r[X(n)]\leq(3(1-p_n)^2)^{rn}\{1+3^{-n}2rn(1-p_n)^{-r^2}\}
\end{align}
Since $r$ is fixed and $\lim_{n\rightarrow
\infty}(3(1-p_n)^2)^n=\lambda $, then:
\begin{align}
\lim_{n\rightarrow \infty}(E_r[X(n)])=\lambda^r \text{ for
r=0,1,2,...}
\end{align}
This shows that $ X(n)\stackrel{d}{\longrightarrow}\textbf{P}_\lambda$. \done\\\\

\noindent{\textbf{For $p=0.5$:}}\\
In the calculation of $p_n=1-1/\sqrt3\lambda^{1/2n}(1+o(1/n))$, if
we fix $p=(\sqrt3-1)/(\sqrt3)$ and let $\lambda=1$, then we get
that the distribution of $X(n)$ , as $n\rightarrow \infty$, tends
to a Poisson distribution with mean $1$. Therefore, one can
conclude:
\begin{align}
\lim_{n\rightarrow \infty}(1-g_n(p))=\lim_{n\rightarrow
\infty}(\textbf{P}(X(n)=0))=e^{-1}
\end{align}
\noindent This shows that for $p=(\sqrt3-1)/(\sqrt3)$ the
probability that $^3Q_{p_n}^n$ has no isolated point, as $n\rightarrow
\infty$, tends to $e^{-1}$. \done

\section{Isoperimetric problem for the 3-cube}
One needs an inequality which relates the size of a subgraph of
$^3Q^n$ to the size of its boundary. This inequality will be
applied to prove that for $p\geq (\sqrt3-1)/(\sqrt3)$, the
probability that subgraphs of $^3Q^n$ do not have a component of
size larger than $2$ and smaller than $3^n$, as $n\rightarrow \infty$, tends to $1$.

\noindent\begin{definition}The edge boundary $b_G(H)$, where $H$
is an induced subgraph of G, is the number of edges which joins
vertices in $H$ to the vertices in $G\backslash
H$.\end{definition}

\noindent The main task in this section is to calculate
$b_{^3Q^n}(k)$. Since $^3Q^n$ is 2n-regular and $H$ is an induced
subgraph of $G$ with $|V(H)|=m$, then:
\begin{align}
b_{^3Q^n}(H)& = mn-2e(H),   \text{where e(H) is the total number of edges in H.}\\
b_{^3Q^n}(k)& = mn-2e_n(k), \text{where }e_n(k)=max\{ e(H):H
\text{ induced subgraph of } ^3Q^n, |V(H)|=m \}.\end{align}
\begin{definition}
$h(i):=$  sum of digits in the ternary expansion of $i$ and
$f(l,m):=\displaystyle\sum_{l\leq i<m}h(i)$
\end{definition}

\noindent\begin{lemma}\label{2.1} If $1 \leq k \leq l$, then
$f(l,l+k)\geq f(0,k)+k$
\end{lemma}

\proof \noindent Look at the ternary expansion of a few numbers:
\[
\begin{array}{ccccc}
Column&2&1&0\\
Bin\backslash Dec&3^2&3^1&3^0\\
  0    &   &   & 0  \\
  1    &   &   & 1  \\
  2    &   &   & 2  \\
  3    &   & 1 & 0  \\
  4    &   & 1 & 1 \\
  5    &   & 1 & 2 \\
  6    &   & 2 & 0 \\
  7    &   & 2 & 1 \\
  8    &   & 2 & 2 \\
  9    & 1 & 0 & 0 \\
 10    & 1 & 0 & 1 \\
 11    & 1 & 0 & 2\\
 12    & 1 & 1 & 0   \\
 13    & 1 & 1 & 1  \\
 14    & 1 & 1 & 2 \\
 15    & 1 & 2 & 0 \\
 16    & 1 & 2 & 1 \\
 17    & 1 & 2 & 2 \\
 18    & 2 & 0 & 0 \\
 19    & 2 & 0 & 1 \\
 20    & 2 & 0 & 2 \\
 21    & 2 & 1 & 0
\end{array}
\]

From this representation, one can observe that column $i$ starts
with a block of $3^i$ zeros. Therefore, sum of jth digits of $k$
numbers, one after the other, is minimal if the first block of
$0$'s is as long as possible. Hence, one can conclude:
\begin{align}\label{2.2}
f(l,l+k)\geq f(0,k)
\end{align}

For every $i$ define $r$ such that $ 0\leq i \leq 3^r-1$. The
ternary expansions of $i$ and $3^r-1-i$ are symmetric. This means
that, if there is a $2/0$ in an specific location of the ternary
expansion of $i$ then there is a $0/2$ in the same location of the
ternary expansion of $3^r-1-i$; and if there is a $1$ in an
specific location of the ternary expansion of $i$ then there is a
$1$ in the same location of the ternary expansion of $3^r-1-i$.
Therefore,
\begin{align}
h(i)+h(3^r-1-i)=2r \text{ for } 0\leq i \leq 3^r-1
\end{align}
\noindent Consequently, since:
\begin{align}
\displaystyle\sum_{l\leq i<l+k}h(i)+\displaystyle\sum_{3^r-l-k\leq
i <3^r-l}h(i)=2rk
\end{align}
\noindent then:
\begin{align}\label{2.5}
f(l,l+k)+f(3^r-l-k,3^r-l)=2rk, \text{where } l+k\leq 3^r
\end{align}

\noindent One should first prove lemma \ref{2.1} with the
assumption $k\leq 3^r \leq l$ by using \ref{2.2} and \ref{ 2.5}.
This assumption means that the length of the sequence of 0's, 1's
and 2's in the ternary expansion of $3^r+k$ and $3^r$ are equal.\\

\noindent With the same logic that one gets \ref{2.2}, one gets:
\begin{align}
f(l,l+k)\geq f(3^r,3^r+k) \text{, when } 3^r\leq l
\end{align}
\noindent and then for $k \leq 3^r$ one can get:

\begin{align} \label{2.3}
f(3^r,3^r+k)=\displaystyle\sum_{3^r\leq i < 3^r+k}h(i)
\end{align}

\noindent $\sum_{3^r\leq i<3^r+k}h(i)$ is the sum over numbers
with the same length in their ternary expansion's sequence. When
one removes the last digit in their ternary expansion, what
remains is $f(0,k)$. Therefore, $\sum_{3^r\leq i<3^r+k}h(i)$ is
equal to sum of the last digits plus $f(0,k)$. Hence:
\begin{align} \label{2.4}
f(l,l+k)\geq f(3^r,3^r+k)=k+f(0,k) \text{ when } k\leq 3^r \leq l
\end{align}
and finally:
\begin{align} \label{2.1.24}
f(l,l+k)\geq f(0,k)+k \text{ where }k\leq3^r\leq l
\end{align}
\noindent Lemma \ref{2.1} is proved with the assumption $k\leq 3^r
\leq l$. Now, one should prove lemma \ref{2.1} without this
assumption. The proof is based on induction on $K$. We want to prove that for
$1\leq K\leq l$, $f(l+K,l) \geq K + f(0,K)$. Fix $k$ such that $1\leq k \leq l$ and $K<k$. For $K=1$ the
inequality is trivial. Assume that the inequality is true for
$K<k$ and $K>2$, which means:
\begin{align} \label{2.7}
f(l,l+K)\geq K+f(0,K) \text{ when } 1\leq k \leq l \text{ , and }
K<k
\end{align}
Now, one should prove the inequality for $K=k$. Define $r\geq1$ by
$3^{r-1}\leq k<3^r$. If $l\geq3^{r}$, then $k\leq 3^r \leq l$ and the lemma is implied by
inequality \ref{2.1.24}. One may assume that $3^{r-1}< l < 3^r$.
Finally, one should apply inequality \ref{2.5} and \ref{2.7} in
order to get the final result:
\begin{align}
f(l+k)&=f(l,3^r)+f(3^r,l+k) \text{\emph{ (from definition of f and }} l\geq 3^r) \\
&=(3^r-l)2r-f(0,3^r-l)+f(3^r,l+k) \text{ \emph{ (from \ref{2.5})}} \\
&\geq (3^r-l)2r-f(0,3^r-l) +f(0,l+k-3^r)+l+k-3^r \text{\emph{ (from \ref{2.7})}} \\
&\geq (3^r-l)2r-f(3^r-k,3^r-k+3^r-l)+3^r-l +f(0,l+k-3^r)+l+k-3^r \text{ \emph{(from \ref{2.7})}} \\
&\geq (3^r-l)2r-f(3^r-k,3^r-k+3^r-l)+f(0,l+k-3^r)+k  \\
&\geq f(l+k-3^r,k)+f(0,l+k-3^r)+k \text{ \emph{(from \ref{2.5})}} \\
&\geq f(0,k)+k \text{ \emph{(from characteristics of f)}}
\end{align}
\done
\begin{theorem} \label{t.2.2}
For $2\leq m \leq 3^{n}$ we have $b_{^3Q^n}(m)=mn-2f(0,m)$. In
other words, $f(0,m)=e_n(m) \text{ where }e_n(m)=\max\{ e(H):H
\text{ induced subgraph of $^3Q^n$ } |V(H)|=m \}.$
\end{theorem}
\proof First, let us fix an $m$. As the first step, one should
prove that $e_{n}(m)\geq f(0,m)$. Vertex $i$ is connected to
$h(i)$ vertices $j$ with $j<i$, since for each $1$ ($2$) in the
ternary expansion of $i$ there is exactly one j $(j<i)$, which its
ternary expansion differs in the position of that $1$ ($2$).
Therefore, on can conclude that $W=\{0,1,2,...,m-1 \}$ contains
$\sum_{0\leq i<m}h(i)=f(0,m)$ edges. So, $e_{n}(m)\geq f(0,m)$.

As the second step, one should prove that $e_n(m)\leq f(0,m)$ by
induction on $n$ for fixed $m$. For $n=1$ the inequality is
trivially true. Assume that it is true for $N<n$, which means:
\begin{align}
e_N(m)\leq f(0,m), \text{ where } N<n \text{ and the fixed m is }
3\leq m \leq 3^{n}\label{2.21}
\end{align}

\noindent Now, one should check the inequality \ref{2.21} for
$N=n$. This means that we should find an $H$ induced subgraph of
$^3Q^n$, $|V(H)|=m$, which maximize $e_n(m)$. Let us split $^3Q^n$
into three (n-1)-dimensional cubes, face-1, face-2 and face-3 each
with $3^{n-1}$ vertices and $(n-1)3^{n-1}$ edges. Now, one can
construct $H$. Choose $m_1$ vertices for $H$ from face-1, $m_2$
vertices from face-2, and $m_3$ vertices from face-3 where
$m_1+m_2+m_3=m$ and $m_1\leq m_2\leq m_3$ and $m_1+m_2\leq m_3$.
In other words, $H$ is constructed from three induced subgraphs,
denoted by $H_1$, $H_2$, and $H_3$.

Each face is a (n-1)-dimensional cube, therefore inequality
\ref{2.21} holds for $H_1$, $H_2$ and $H_3$. Also, each vertex of
each face is connected to exactly on vertex from one face and one
vertex from the other face. Hence, the number of vertices of $H$
is at most:
\begin{align}
e_{n}(m) \leq f(0,m_1) + f(0,m_2) + f(0,m_3) + 2m_1 + m_2
\end{align}
\noindent $2m_1+m_2$ is the maximum number of edges between three
faces that one can choose here. $2m_1$ is the maximum number of
edges between chosen vertices in $H_1$ and $H_2$ plus the maximum
number of edges between chosen vertices in $H_1$ and $H_3$.
Consequently, $m_2$ is the maximum number of edges between chosen
vertices in $H_2$ and $H_3$. Therefore:
\begin{align}
e_{n}(m) & \leq f(0,m_1) + f(0,m_2) + f(0,m_3) + 2m_1 + m_2 \text{ (from \ref{2.21})}\\
& \leq f(m_2,m_2+m_1)+f(0,m_2) +m_1+m_2\text{ (from lemma \ref{2.1})} \\
& \leq f(m_3,m_3+m_2+m_1)+f(0,m_3)\text{ (from lemma \ref{2.1})} \\
& \leq f(0,m) \text{ (from definition of f)}
\end{align}
\done\\
 Theorem \ref{t.2.2} shows that, if we want to choose an
induced subgraph of $^3Q^n$, with $m$ vertices, which has the
smallest edge boundary, then we should choose the induced subgraph
of $^3Q^n$ with the set of vertices $W=\{0,1,2,...,m-1\}$.\\

\begin{corollary}\label{2.6}
For all $k$ and $n$, $e_n(k)\leq k\lceil log_3k\rceil$, which is
equivalent to $b_{^3Q^n}(k)\geq 2k(n-\lceil log_3k\rceil)$.
\end{corollary}
\proof \noindent Let $r=\lceil \log_3k\rceil$. Then
\begin{align}
2f(0,k)&=f(0,k)+f(0,k)\leq f(0,k)+f(3^r-k,3^r) \text{ (from
\ref{2.2})}\\
&=2rk \text{ (from \ref{2.5})}
\end{align}
\noindent Therefore:
\begin{align}
e_{n}(k)=f(0,k)\leq rk=k\lceil \log_3k \rceil
\end{align}
\done

\section{Isolated components of size larger than 2 and smaller than $3^n$}
\begin{definition}
$C_s$ is the family of s-subsets (subsets with size s) of
$V=V(^3Q^n)$ whose induced graph is connected.
\end{definition}

\noindent\textbf{Remarks:} $h(n):=o(g(n))$ means $\frac{h(n)}{g(n)}\rightarrow 0$ as $n\rightarrow \infty$.\\
\noindent\textbf{Remarks:} The following inequality will be
applied a lot in the rest of this section:
\begin{align}
(\frac{n}{k})^k \leq \binom{n}{k} \leq \frac{n^k}{k!} \leq
(\frac{ne}{k})^k \label{rem.2.1}
\end{align}

\begin{theorem}
If $p_n\geq 1- \frac{1}{\sqrt 3}(\log n)^{\frac{1}{n}}$, the
probability that for some $S\in C_s, 2\leq s \leq \frac{3^n}{2}$,
no edges of $^3Q_{p_n}^n$ join $S$ to $V(^3Q^n)\setminus S$, as $n
\rightarrow \infty$, tends to 0.
\end{theorem}
\noindent \textbf{Note:} For $\frac{3^n}{2}<s<3^n$, if there exist
a component of size smaller than $3^n$ then there is at least one
component of size smaller than $\frac{3^n}{2}$ which contradicts
with the theorem.
 \proof \noindent Consider $S\subset V=V(^3Q^n)$ and
set $b(S)=b_{Q^n}(H)$  where $H$ is the induced subgraph of
$^3Q^n$ with the set of vertices S. One can observe that:
\begin{align}
\textbf{P}(\text{No edges of } ^3Q_{p_n}^n \text{ join S to }
V\setminus S)= (1-p_n)^{b(S)}
\end{align}
\noindent In order to prove the theorem, it is sufficient to show:
\begin{align}
\sum_{s=2}^{3^{\frac{n}{2}}}\sum_{S\in C_s}(1-p_n)^{b(S)}=o(1)
\end{align}
From corollary \ref{2.6}, one knows that for $|S|=s$:
\begin{align}
b(S)\geq b(s)\geq 2s(n-\lceil \log_3s \rceil)
\end{align}
\noindent and therefore:
\begin{align}
\sum_{S\in C_s}(1-p_n)^{b(S)}\leq |C_s|(1-p_n)^{b(s)}
\end{align}
\noindent One may partition $s$, $2\leq s \leq \frac{3^n}{2}$, to
different  intervals in order to find a small enough bound for $|C_s|$ and $(1-p_n)^{b(s)}$.\\\\
\textbf{First interval $2 \leq s \leq s_1, s_1= \lfloor
\frac{3^{\frac{n}{2}}}{n^2} \rfloor$:}\\

First, one should find a bound for $|C_s|$. One has maximum $3^n$
choices to choose the first element for $C_s$. The selected element
is connected to maximum $2n$ vertices, therefore, there are
maximum $2n$ choices to choose the second element. With the same
logic, there are at most $2n(s-1)$ choices to choose the last element for
$|C_s|$. Hence, one can show:
\begin{align}
|C_s|\leq 3^n(2n)(2n(2))...(2n(s-1)) \leq (s-1)!(2n)^{s-1}3^n
\end{align}
\noindent and:
\begin{align}
|C_s|(1-p_n)^{b(s)}\leq (s-1)!(2n)^{s-1}3^n (1-p_n)^{2s(n-\lceil
\log_3s \rceil)} \label{2.eq.2}
\end{align}
Since $p_n=1-\frac{1}{\sqrt3}(\log n)^{\frac{1}{n}}$, so for large enough $n$:
\begin{align}
(1-p_n)^{2s(n-\lceil log_3s \rceil)}&\leq(3)^{-ns}(\log n)^{2s}(1-p_n)^{-2s(\log_3s)} \text{(neglecting some small terms)}\\
&\text{( since: }(\log n)^{\frac{-2s\log_3 s}{n}} \leq 1 \text{ for large enough }n) \\
& = (3)^{-ns}(\log n)^{2s}3^{s\log_3s}(\log n)^{\frac{-2s\log_3 s}{n}}\\
 &\leq (3)^{-ns}(\log n)^{2s}s^s \label{2.eq.1}
\end{align}
From equations \ref{2.eq.2} and \ref{2.eq.1}, one gets:
\begin{align}
|C_s|(1-p_n)^{b(s)}\leq (s-1)!(2n)^{s-1}3^n(3)^{-ns}(\log
n)^{2s}s^s \label{2.eq.3}
\end{align}
\noindent Assume that the right hand sides of inequality
\ref{2.eq.3} is equal to A. After multiplying both side of
inequality \ref{2.eq.3} with $\frac{2ns^{s+1}}{s!}$ and then
getting $\log_3$ from both sides, one gets:
\begin{align}
\log_3(|C_s|(1-p_n)^{b(s)}\frac{2ns^{s+1}}{s!})\leq\log_3(A\frac{2ns^{s+1}}{s!})
\end{align}
\noindent If $\log_3(A\frac{2ns^{s+1}}{s!})\rightarrow-\infty$ as
$n\rightarrow\infty$ then $A\frac{2ns^{s+1}}{s!}$ should tend to
0. This means that $|C_s|(1-p_n)^{b(s)}\frac{2ns^{s+1}}{s!}$ tends
to 0, as $n\rightarrow \infty$. Therefore:
\begin{align}
|C_s|(1-p_n)^{b(s)}\leq \frac{s!}{2ns^{s+1}} \text{ for large
values of n}
\end{align}
which shows that:
\begin{align}
\sum_{s=2}^{s_1}\sum_{S\in C_s}(1-p_n)^{b(S)}=o(1)
\end{align}

\noindent Finally, it remains to prove
$\log_3(A\frac{2ns^{s+1}}{s!})\rightarrow-\infty$ as
$n\rightarrow\infty$. One can verify this for $s\leq n$ and $s>n$.\done \\\\

\noindent \textbf{Second interval $s_1 +1 \leq s \leq
\frac{3^n}{2}$ and
$S\in C_s^-, s_1= \lfloor \frac{3^{\frac{n}{2}}}{n^2} \rfloor$:}\\
Define $C_s^-$ and $C_s^+$ as follows:
\begin{definition}
\begin{align}
C_s^-:=\{S\in C_s | b(s) \geq 2s(n-\log_3s+\log_3n)\} \text{, and
} C_s^+:=C_s\backslash C_s^-
\end{align}
\end{definition}
One can bound $|C_s^-|$ for  $s_1 +1 \leq s \leq 3^\frac{n}{2}$  as
follows:
\begin{align}
|C_s^-|\leq |C_s| \leq \binom{3^n}{s} \leq \frac{3^{ns}}{s!} \leq
(\frac{e3^n}{s})^s
\end{align}
Hence:
\begin{align}
\displaystyle\sum_{s=s_1+1}^{3^{\frac{n}{2}}}\sum_{S\in C_s^-}(1-p_n)^{b(S)} &\leq \sum_{s=s_1+1}^{3^{\frac{n}{2}}}(\frac{e3^n}{s})^s(\sqrt3^{-1}(\log n)^{\frac{1}{n}})^{2s(n-\log_3s+\log_3n)}\\
\displaystyle& \leq \sum_{s=s_1+1}^{3^{\frac{n}{2}}}(\frac{e3^n3^{-(n-\log_3 s +\log_3 n)}(\log n)^{\frac{2(n-\log_3s+\log_3n)}{n}}}{s})^s\\
& \leq \sum_{s=s_1+1}^{3^{\frac{n}{2}}}(\frac{e3^n3^{-n}3^{\log_3s}3^{-\log_3n}(\log n)^2}{s})^s(\log n)^{\frac{2n(-\log_3s+\log_3n)}{s}}\\
& \text{ (since for large enough n: } (\log n)^{\frac{2n(-\log_3s+\log_3n)}{s}} \leq1 )\\
& \leq \sum_{s=s_1+1}^{3^{\frac{n}{2}}} (\frac{e(\log n)^2}{n})^s = o(1)
\end{align}
\done\\\\
\noindent \textbf{Third interval $s_1 \leq s \leq s_2,s_1= \lfloor
\frac{3^{\frac{n}{2}}}{n^2} \rfloor, s_2=\lfloor \frac{3^n}{(\log n)^7}\rfloor$ and $S\in C_s^+$:}\\
For the 3rd and the 4th intervals one needs to find a bound for
$|C_s^+|$. The following lemma, presented by B.Bollobas
\cite{bol}, helps us in this matter:
\begin{lemma}\label{2.l.1}
Let G be a graph of order $v$ and suppose that $\Delta(G)\leq
\Delta$, $2e(G)=vd$ and $\Delta +1 \leq u \leq v- \Delta -1$.
Then, there is a u-set of U of vertices with \cite{bol}:
\begin{align}
|N(U)|=|U \cup \Gamma (U)| \geq v\frac{d}{\Delta} \{
1-\exp(\frac{-u(\Delta+1)}{v}) \}
\end{align}
where, $\Delta(G):=$ Maximum degree in G, $d:=$ average degree in
G and $\Gamma(U)=\{x\in V(G): xy \in E(G) \text{ for some y}\in
U\}$
\end{lemma}
Let $H=$$^3Q^n[S]$ (the induced subgraph of $^3Q_n$ with the set
of vertices $S$). From the definition of $C_s^+$ one knows that
the average degree in $H$ is at least:
\begin{align}
2(\log_3s-\log_3n)
\end{align}
The goal is to find $U\subset S$, where
$|U|:=u:=\lfloor\frac{2s}{n}\rfloor$, $\Delta = 2n$, $v=s$, $d\geq
\log_3s-\log_3n$ and then use lemma \ref{2.l.1} to calculate the
boundary size of $U,|N(U)|$. First, check if $2n+1 \leq
\lfloor\frac{2s}{n}\rfloor$, as $n\rightarrow\infty$:
\begin{align}
&\frac{2s}{n}=\frac{23^{n/2}}{n^3} \text{for minimum s, and trivially }  2n+1\leq \frac{23^{n/2}}{n^3}, \text{for large enough }n\\
\end{align}
\noindent and then check if $\lfloor\frac{2s}{n}\rfloor \leq
s-(2n+1)$. One should check if  $ns-n(2n+1)\geq 2s$ which means
one should check that whether:
\begin{align}\frac{3^{\frac{n}{2}}(n-2)}{n^3(2n+1)}\geq1\end{align}
\noindent which is clearly true for large enough $n$. Now, one can
apply lemma \ref{2.l.1} on the graph generated by $S$ and get:
\begin{align}
\exists U \subset S: |N(U)|\geq
s\frac{2(\log_3s-\log_3n)}{2n}\{1-\exp(-\frac{u(2n+1)}{s})\}
\end{align}

\noindent where:
\begin{align}
\frac{2(\log_3s-\log_3n)}{2n}\geq \frac{(\log_3(\frac{3^{\frac{n}{2}}}{n^2})-\log_3n)}{n}=\frac{\frac{n}{2}-3\log_3n}{n}\\
\text{which }
\lim_{n\rightarrow\infty}\frac{\frac{n}{2}-3\log_3n}{n}=
\frac{1}{2} \label{2.8}
\end{align}
\noindent on the other hand:
\begin{align}
\lim_{n\rightarrow
\infty}(1-\exp(-\frac{2n+1}{s}(\frac{2s}{n}+1)))=1-e^{-4}
\label{2.9}
\end{align}
Therefore, from \ref{2.8} and \ref{2.9} one gets:
\begin{align}
|N(U)| \geq \frac{1}{2}(1-e^{-4})s \geq \frac{s}{3} \text{ as }
n\rightarrow \infty\label{2.10}
\end{align}

\noindent Now that we have $|N(U)|$, we can estimate a bound for
$|C_s^+|$ here. We know from \ref{2.10} that for each $S \in
C_s^+$ there exist a $U \subseteq S$,
$|U|:=u:=\lfloor\frac{2s}{n}\rfloor$, such that
$|N(U)|\geq s/3$. Therefore, one can choose $S\in C_s^+$ as follows:\\
\noindent 1. Select u vertices of $^3Q^n$; there are $\binom{3^n}{u}$ choices for this u.\\
\noindent 2. Select $\lfloor\frac{s}{3}\rfloor-u$ neighbors of the
selected vertices of u in part 1; there are maximum $(2^{n})^u$
choices, since there are at most
$\binom{2n}{0}+\binom{2n}{1}+\binom{2n}{2}+...\binom{2n}{2n}=2^{2n}$
ways to
find neighbors of a vertex in $U$.\\
\noindent 3. Select $\lfloor\frac{2s}{3}\rfloor$ other vertices; there are at most $\binom{3^n}{\lfloor\frac{2s}{3}\rfloor}$ choices.\\
\noindent Hence:
\begin{align}
|C_s^+|\leq
\binom{3^n}{u}(2^{2n})^u\binom{3^n}{\lfloor\frac{2s}{3}\rfloor}
\end{align}
\noindent \noindent and:
\begin{align}
\sum_{S\in C_s^+}(1-p_n)^{b(S)}\leq
\binom{3^n}{u}(2^{2n})^u\binom{3^n}{\lfloor\frac{2s}{3}\rfloor}(1-p_n)^{b(s)}\label{2.11}
\end{align}
\noindent where:
\begin{align}
(1-p_n)^{b(s)} \leq 3^{-sn}s^s(\log n)^{2s} \label{2.12}
\end{align}
\noindent consequently from \ref{2.11}, \ref{2.12} and
\ref{rem.2.1}:
\begin{align}
\sum_{S\in C_s^+}(1-p_n)^{b(S)}\leq
(\frac{e3^n}{u})^u2^{2un}(\frac{e3^n}{\lfloor\frac{2s}{3}\rfloor})^{\lfloor\frac{2s}{3}\rfloor}3^{-sn}s^s(\log
n)^{2s}\label{2.14}
\end{align}
\noindent Write $s=3^{\beta n}$, ($\beta = \frac{log_3s}{n}$), so
that:
\begin{align}
3^{\beta n} \leq \frac{3^n}{(\log n)^7}  \Rightarrow \beta \leq
1-\frac{7\log_3\log n}{n} \label{2.13}
\end{align}
\noindent Now, find a bound for the inequality \ref{2.14}. First
calculate the first part of the inequality:
\begin{align}
(\frac{e3^n}{u})^u 2^{2un}
(\frac{e3^n}{\lfloor\frac{2s}{3}\rfloor})^{\lfloor\frac{2s}{3}\rfloor}&
\leq (\frac{e3^n}{\frac{2s}{n}})^{\frac{2s}{n}} 3^{2s} 2^{4s}
(\frac{e3^n}{\frac{2s}{3}})^{\frac{2s}{3}}\\
&\text{ (since for large enough n and $s_1 \leq s \leq s_2$: }(\frac{e3^n}{\frac{2s}{n}})^{\frac{2s}{n}}\leq1)\\
&\leq(3^22^4(\frac{3}{2}e)^{\frac{2}{3}})^s\frac{3^{\frac{2s}{3n}}}{s^{\frac{2s}{3}}}
=c^s3^{\frac{2}{3}sn(1-\frac{\log_3s}{n})}=c^s3^{\frac{2}{3}sn(1-\beta)}\label{2.15}
\end{align}
\noindent where c is a positive constant. Now, by substituting
\ref{2.15} in \ref{2.14} one gets:
\begin{align}
\sum_{S\in C_s^+}(1-p_n)^{b(S)} & \leq 3^{-sn}s^s(\log n)^{2s}c^s3^{\frac{2}{3}sn(1-\beta)} \\
&=c^s(\log n)^{2s}3^{-\frac{sn(1-\beta)}{3}}\\
&\leq c^s(\log n)^{2s}3^{-\frac{7s\log_3\log n}{3n}} \text{ , (from \ref{2.13})} \\
&=c^s(\log n)^{2s}3^{\log_3(\log n)^{\frac{-7s}{3}}}\\
&\leq c^s(\log n)^{2s}(\log n)^{\frac{-7s}{3}}\\
& =c^s(\log n)^{\frac{-s}{3}} \label{2.28}
\end{align}
\noindent and finally from \ref{2.28}:
\begin{align}
\sum_{s=s_1}^{s_2}\sum_{S\in C_s^+}(1-p_n)^{b(S)}\leq
\sum_{s=s_1}^{s_2} c^s(\log n)^{\frac{-s}{3}} = o(1)
\end{align}
\done\\\\
\noindent \textbf{Fourth interval $s_2+1 \leq s \leq \frac{3^{n}}{2}$ and $,s_2=\lfloor \frac{3^n}{(\log n)^9}\rfloor, S\in C_s^+$:}\\

In $H=$$^3Q^n[S]$ (the induced subgraph $^3Q^n$ with the set of
vertices S), the average degree is at least:
\begin{align}
2(\log_3s-\log_3n)>2(n-2\log_3n)\label{2.33.1}
\end{align}
\noindent since:
\begin{align}
s\geq \lceil\frac{3^n}{(\log n)^9}\rceil \Rightarrow
\log_3(\frac{3^n}{(\log n)^9})<\log_3 s\\
\Rightarrow \log_3s - \log_3 n \geq \log_3(\frac{3^n}{(\log
n)^9})-\log_3n \geq n -\log_3(\log n)^9-\log_3n\\
(\text{for large enough n one can get } n > (\log n)^9)\\
\geq n-2\log_3n
\end{align}
\noindent First, look for a subgraph of H with large average
degree. Let T be the set of vertices of H with degree at least
$2(n-(\log_3n)^2)$, and set $t=|T|$. From \ref{2.33.1} one can
conclude that the sum of degrees in $H$ is at least
$s(n-2\log_3n)$. We also know that:
\begin{align}
\text{Sum of degrees in }S & \leq 2s(n-2\log_3n)\\
&\leq t\times(\text{Maximum degree of vertices in set $T$ of graph $H$ })\\
&+ (s-t)\times(\text{Maximum degree of vertices in set $S\setminus
T$ of
graph $H$ })\\
& \leq 2tn +2(s-t)(n-(\log_3n)^2\\
\Rightarrow t\geq s(1-\frac{4}{\log_3n})\label{2.34.1}
\end{align}
\noindent Define $H_1=$$^3Q^n[T]=H[T]$ as the induced subgraph
spanned by $T$. We want to calculate $|N_{H_1}(U)|$ in $H_1$,
hence we should estimate the size of $H_1$ and after that
calculate the average degree in T. Let us first calculate
$e(H_1)$, the total number of edges in $H_1$.
\begin{align}
e(H_1)\geq e(H)-2(s-t)n
\geq2\frac{s}{2}(n-2\log_3n)-\frac{4s}{\log_3n}n \text{ (from
\ref{2.33.1} and \ref{2.34.1})}
\end{align}
\noindent One knows that the average degree in $H_1$ is at least
$\frac{2e(H_1)}{s}$, and:
\begin{align}
\frac{2e(H_1)}{s}\geq2(n-2\log_3n)-\frac{8n}{\log_3n}\geq
2n-\frac{9}{\log_3n}\\
\text{(since: } \log_3n^2<\frac{n}{\log_3n} \text{ for large
enough n)}
\end{align}
\noindent Set $u=\lfloor\frac{3^n}{n^{\frac{1}{2}}}\rfloor$. One
should check the conditions of lemma \ref{2.1} here. Let
$v=t,\Delta=2n, d\geq 2n-\frac{9}{\log_3n}$. So, one should check
if $2n+1\leq\frac{3^n}{n^{\frac{1}{2}}}\leq t-2(n+1)$, for large
enough $n$. Clearly, $2n+1\leq\frac{3^n}{n^{\frac{1}{2}}}$, as
$n\rightarrow \infty$. It remains to prove
$\frac{3^n}{n^{\frac{1}{2}}}\leq t-2(n+1)$, for large enough $n$.
For minimum $s$ from \ref{2.34.1} we can get:
\begin{align}
t\geq& \frac{3^n}{(\log n)^7}(1-\frac{4}{\log_3n}) \text{ (from \ref{2.34.1})}\\
&\geq \frac{3^n}{n^\frac{1}{2}} +2(n+1)\text{ (for large enough n)
}\end{align}
 Now, one can use lemma \ref{2.1} and estimate $|N_{H_1}(U)|$.
\begin{align}
|N_{H_1}(U)|&\geq
    \frac{t}{2n}(2n-\frac{9n}{\log_3n})\{1-\exp(-\frac{2n+1}{t}\frac{3^n}{n^{\frac{1}{2}}})\}\\
&\geq\frac{t}{2}(n-\frac{9}{\log_3n})\{1-\exp(-\frac{2n+1}{t}\frac{3^n}{n^\frac{1}{2}})\}
\label{2.18}
\end{align}
\noindent Let us estimate a bound for
$\exp(-\frac{2n+1}{t}\frac{3^n}{n^{\frac{1}{2}}})$. One knows that
$t\geq s(1-\frac{4}{\log_3n})$. Since $\max(t)=s$ and
$\max(s)=\frac{3^{n}}{2}$, then:
\begin{align}
\frac{3^n(2n+1)}{n^{\frac{1}{2}}t}\geq
\frac{3^n(2n+1)}{n^{\frac{1}{2}}3^{n}} =
\frac{(2n+1)}{n^{\frac{1}{2}}}\geq n^\frac{1}{4}\text{ (for large
enough n)
}\\
\Rightarrow \{1-\exp(-\frac{2n+1}{t}\frac{3^n}{n^\frac{1}{2}})\}
\geq \exp(-n^\frac{1}{4})\text{ (for large enough n) }\label{2.19}
\end{align}
\noindent By using the bound from \ref{2.19} in \ref{2.18}, one
gets:
\begin{align}
|N_{H}(U)|\geq|N(H_1)|& \geq \frac{t}{2}(2-\frac{9}{\log_3n})\{
1-\exp(-n^{\frac{1}{4}})\}\\
&=\frac{t}{2}\{2+\exp(-n^{\frac{1}{4}})\frac{9}{\log_3n}-2\exp(-n^{\frac{1}{4}})-\frac{9}{\log_3n}\}\\
&(\lim_{n\rightarrow\infty}\exp(-n^{\frac{1}{4}})\frac{9}{\log_3n}=0
\text{ and } \exp(-n^{\frac{1}{4}})< \frac{1}{\log_3n} \text{ (for
large enough n)
})\\
&\geq
\frac{t}{2}\{2-\frac{2}{\log_3n}-\frac{9}{\log_3n}\}=\frac{t}{2}(2-\frac{11}{\log_3n})\\
&\geq\frac{s}{2}(1-\frac{2}{\log_3n})(2-\frac{11}{\log_3n})=\frac{s}{2}(2+\frac{2}{\log_3n}\frac{11}{\log_3n}-\frac{15}{\log_3n})
\text{ (for large enough n)
}\\
&\geq\frac{s}{2}(2-\frac{15}{\log_3n})=s(1-\frac{7.5}{\log_3n})\label{2.20}
\end{align}
\noindent Now that we have $|N_{H}(U)|$, we can estimate a bound
for $|C_s^+|$ here. We know from \ref{2.20} that for each $S \in
C_s^+$ there exist a $U \subseteq S$,
$|U|:=u:=\lfloor\frac{3^n}{n^\frac{1}{2}}\rfloor$, such that
$|N_{H}(U)|\geq s(1-\frac{7.5}{\log_3n})$. Therefore, one can choose $S\in C_s^+$ as follows:\\
\noindent 1. Select u vertices of $^3Q^n$; there are $\binom{3^n}{u}$ choices for this u.\\
\noindent 2. Select $\lfloor s(1-\frac{7.5}{\log_n3})\rfloor-u$
neighbors of the selected vertices in part 1. At most
$2(\log_3n)^2$ of the $2n$ neighbors of a vertex in $U$ do not
belong to $N_{H}(U)$. Hence there are at most
 $\sum_{(k_j)}(\prod_{i=1}^{u}\binom{2n}{j})$ ways to find
neighbors of $u$ vertices in $U$, where the sum is over all
$(k_1,k_2,...,k_u), k_i \leq 2(\log_3n)^2$. We know that:
\begin{align}
\sum_{(k_i)}\prod_{i=1}^{u}\binom{2n}{k_i}  &\leq \sum_{(k_i)}\prod_{i=1}^{u}(\frac{(2n)^i}{i!})\\
&\leq \sum_{(k_i)}\prod_{i=1}^{u}(\frac{(2n)^{(2(\log_3n)^2)}}{(2(\log_3n)^2)!})\\
&= \sum_{(k_i)}\frac{(2n)^{2u(\log_3n)^2}}{((2(\log_3n)^2)!)^u}\\
&= {(2(\log_3n)^2)}^u\frac{(2n)^{2u(\log_3n)^2}}{((2(\log_3n)^2)!)^u}\\
& \leq (2n)^{2u(\log_3n)^2}\\
\end{align}

\noindent 3. Select $\lfloor\frac{7.5}{\log_3n}\rfloor$ other vertices; there are $\binom{3^n}{\lfloor\frac{7.5}{\log_3n}\rfloor}$ choices.\\
\noindent Hence:
\begin{align}
\sum_{S\in C_s^+}(1-p_n)^{b(S)}\leq
\binom{3^n}{u}(2n)^{2u(\log_3n)^2}\binom{3^n}{\lfloor\frac{7.5}{\log_3n}\rfloor}3^{-2s(n-\log_3s)}(\log
n)^{2s(1-\frac{\log_3s}{n})}\\
\end{align}
\noindent where:
\begin{align}
\binom{3^n}{u}(2n)^{2u(\log_3n)^2}\binom{3^n}{\lfloor\frac{7.5}{\log_3n}\rfloor}
= 3^{o(s)}
\end{align}
\noindent Therefore:
\begin{align}
\sum_{S\in C_s^+}(1-p_n)^{b(S)}\leq 3^{\varepsilon(s)}
\end{align}
\noindent where:
\begin{align}
\varepsilon(s)=o(s)-2s\{n-\log_3s-\log_3\log n
+\frac{\log_3s}{n}\log_3\log n\}
\end{align}
\noindent Since $s\leq \frac{3^{n}}{2}$, hence one can get:
\begin{align}
\varepsilon(s)&\leq o(s)-2s\{n-(n-1)-\log_3\log n
+\frac{n-1}{n}\log_3\log n\}\\
&=o(s)-2s\{1-\frac{1}{n}\log_3\log n\}\leq-s
\end{align}
\noindent Therefore, for large enough n, one can get:
\begin{align}
\sum_{s=s_2+1}^{\frac{3^{n}}{2}}\sum_{S\in
C_s^+}(1-p_n)^{b(S)}\leq
\sum_{s=s_2+1}^{\frac{3^{n}}{2}}3^{-s}=o(1)
\end{align}
\done

\chapter{Connected random subgraph of the $P_3$-product}

The $P_3^n$ graph is the cartesian products of $n$ copies of
$P_3$, where $P_3$ stands for a path with length 2. If one labels
each vertex of $P_3^n$ from $0$ to $3^n-1$, then two vertices are
adjacent if the difference between their ternary representation is
1, in other words, vertex $x=(x_1,x_2,...,x_n)$ is connected to
the vertex $y=(y_1,y_2,...,y_n)$ if for some $i$ we have
$|x_i-y_i|=1$ and $x_j=x_j$ for $j\neq i$. Vertex
$x=(x_1,x_2,...,x_n)$ is connected to $n+i$ vertices where
$i=|\{j|x_j=1, j\in \{1,..,n\}\}|$. A random subgraph of $P_3^n$
contains all vertices of $P_3^n$, and each edge independently with
probability $p$. $p$ is called the percolation parameter and
$P_{3,p}^n$ stands for the random subgraph of $P_3^n$

The main goal in this chapter is to explore a critical value
$p_c$, which for fixed values of $p$ if $p<p_c$ then almost no
random subgraphs of $P_3^n$ is connected, as $n \rightarrow
\infty$; but if $p>p_c$ then almost all random subgraphs of
$P_3^n$ are connected, as $n \rightarrow \infty$. We suggest that
this critical value is $2-\sqrt2$. The proof in this chapter is
not complete and there is a place for further work.

In the first section of this chapter, it is proved that for
$p<2-\sqrt2$ almost no random subgraphs of $P_3^n$ are connected,
as $n \rightarrow \infty$. Then, it is proved that, for
$p>2-\sqrt2$ almost all random subgraphs of $P_3^n$ have no
isolated point, as $n\rightarrow \infty$. In the second section,
the probability that random subgraphs of $P_3^n$, with the
percolation parameter $p>2-\sqrt2$, have no components with size
larger than $1$ and smaller than $3^n$, as $n\rightarrow \infty$,
is explored.

\section{Isolated vertices}
\noindent{\textbf{For $p<2-\sqrt2:$}}\\\\
Let us define $X_i$ and $X$ for a graph $G$ as follows:
\begin{definition}
$X_i(n):=\left\{
\begin{array}{ll}
         1 & \mbox{Vertex $i$ is isolated,  $i\in V(G)$};\\
        0 & \mbox{Vertex $i$ is NOT isolated,  $i\in V(G)$}.\end{array} \right.$
, and $X(n):=\displaystyle\sum_{i\in V(G)}X_i(n)$.
\end{definition}

\noindent As the first step we calculate $E[X(n)]$. One can
categorize the set of vertices $V$ into $n+1$ subsets $V_i$'s
where $x=(x_1,x_2,...,x_n)\in V_i$ if $|\{j|x_j=1, j\in
\{1,..,n\}\}| = i$. Hence:
\begin{align}
\mu := E[X(n)] = \sum_{ j\in V(P^n_{3,p})} E[X_j(n)] = \sum_{
i=0}^n\sum_{ k\in V_i} E[X_k(n)]=
\sum_{i=0}^{n}\binom{n}{i}2^{n-i}(1-p)^{n+i}=(1-p)^{n}(3-p)^n
\label{6.1}
\end{align}
From (\ref{6.1}), the threshold value for $E[X]$ is $p_c=2-\sqrt2$
which is the solution to the equation $(1-p)^n(3-p)^n=1$. This
means that $E[X]\rightarrow \infty$ for $p<p_c$, but
$E[X]\rightarrow 0$ for $p>p_c$. Now, one should calculate
$Var[X]$.
\begin{align}
Var[X(n)] &= \sum_{i\in V(P^n_{3,p})} Var[X_i(n)] + \sum_{i,j\in
V(P^n_{3,p}),i\neq j} Cov[X_i(n),X_j(n)]\\
&= \sum_{ i=0}^n\sum_{ k\in V_i} Var[X_k(n)]+\sum_{i,j\in
V(P^n_{3,p}),i\neq j} Cov[X_i(n),X_j(n)]
\end{align}
\noindent where:
\begin{align}
\sum_{i=0}^n\sum_{k\in V_i}Var[X_k(n)] & = \sum_{i=0}^{n}\binom{n}{i}2^{n-i}((1-p)^{n+i}-(1-p)^{2n+2i}) \\
&=(1-p)^{n}(3-p)^n-(1-p)^{2n}(2+(1-p)^2)^n \label{6.7}
\end{align}
In order to calculate $\sum_{i,j\in V(P^n_{3,p}),i\neq j}
Cov[X_i(n),X_j(n)]$ one should know that if $x\in V_i$ then $x$ is
incident to $2i$ vertices in $V_{i-1}$ and $n-i$ vertices in
$V_{i+1}$. Also, one knows:
\begin{align}
Cov[X_i(n),X_j(n)]& = E[X_i(n)X_j(n)]-E[X_i(n)]E[X_j(n)]=0 \text{
if i,j not adjacent}
\end{align}
Hence:
\begin{align}
\sum_{i,j\in V(P^n_{3,p}),i\neq j}
&Cov[X_i(n),X_j(n)]\\&=\sum_{i=0}^{n}\binom{n}{i}2^{n-i}\{(n-i)(1-p)^{n+i-1}(1-p)^{n+i+1}+2i(1-p)^{n+i-1}(1-p)^{n+i-1}\}
\\
&=\sum_{i=0}^{n}\binom{n}{i}2^{n-i}(1-p)^{2n+2i-1}\{(n-i)(1-p)+\frac{2i}{1-p}\}\\
&=(1-p)^{2n}\sum_{i=0}^{n}\binom{n}{i}2^{n-1}(1-p)^{2i-1}\{n(1-p)+i(\frac{2}{1-p}-(1-p))\}\\
&=n(1-p)^{2n}\sum_{i=0}^{n}\binom{n}{i}2^{n-i}(1-p)^{2i}+(1-p)^{2n-2}(2-(1-p)^2)\sum_{i=0}^{n}\binom{n}{i}2^{n-i}(1-p)^{2i}i\\
&=n(1-p)^{2n}(2+(1-p)^2)^n+(1-p)^{2n-2}(1-(1-p)^2)2^nn\frac{(1-p)^2}{2}(1+\frac{(1-p)^2}{2})^{n-1}\\
&=n(1-p)^{2n}(2+(1-p)^2)^n+n(2-(1-p)^2)(2+(1-p)^2)^{n-1}(1-p)^{2n}\\
&=4n(1-p)^{2n}(2+(1-p)^2)^{n-1} \label{6.8}
\end{align}
 \noindent and finally, from (\ref{6.7}) and (\ref{6.8}) one gets:
\begin{align}
Var[X(n)]&=(1-p)^{n}(3-p)^n-(1-p)^{2n}(2+(1-p)^2)^n+4n(1-p)^{2n}(2+(1-p)^2)^{n-1}\\
&=\mu -
\mu^2\frac{(2+(1-p)^2)^n}{(3-p)^{2n}}+\mu^2\frac{4n(2+(1-p)^2)^n}{(3-p)^{2n}}
\label{6.10}
\end{align}
Now, since we have $Var[X(n)]$, we can use Chebyshev's inequality
to estimate $g_n(p)$. Chebyshev's inequality states that:
\begin{align}
1-\textbf{P}[P_{3,p}^n\text{ contains an isolated
point}]=\textbf{P}[X(n)=0]\leq \textbf{P}[|X(n)-\mu|\geq \mu]\leq
\frac{Var[X(n)]}{\mu^2} \label{6.9}
\end{align}
From \ref{6.10} when $p<2-\sqrt2$, one gets
$Var[X(n)]/\mu^2\rightarrow 0$, as $n\rightarrow\infty$.
Therefore, from (\ref{6.9}) one gets $\lim_{n\rightarrow \infty}
\textbf{P}[P_{3,p}^n\text{ contains an isolated point}]=1$.
Finally, since $$\textbf{P}[P_{3,p}^n\text{ is
connected}]\leq\textbf{P}[P_{3,p}^n\text{ does contains an
isolated point}]$$, then for $p<2-\sqrt2$ the
probability that $P_{3,p}^n$ is connected, as $n\rightarrow \infty$, tends to $0$.\done\\

\noindent{\textbf{For $p>2-\sqrt2$:}}\\\\
In order to calculate $\textbf{P}[P_{3,p}^n\text{ contains an
isolated point}]$ when $p>2-\sqrt2$, as $n\rightarrow \infty$, one
can use the following inequality:
\begin{align}
\textbf{P}[P_{3,p}^n\text{ contains an isolated
point}]=\textbf{P}[X(n)>0]\leq E[X(n)]=\mu
\end{align}
Since $E[X(n)]\rightarrow 0$ as $n\rightarrow\infty$, then
$\lim_{n\rightarrow \infty}\textbf{P}[P_{3,p}^n\text{ contains an
isolated point}]=0$. This means that the probability that there
are no isolated points in $P_{3,p}^n$ for $p>2-\sqrt2$, as
$n\rightarrow
\infty$, tends to $1$.\done\\
\section{Isolated components of size larger than 2 and smaller than $3^n$}
\begin{definition}
$C_s$ is the family of subsets of $V(P_{3}^n)$ with size $s$ whose
their induced graph is connected.
\end{definition}
\noindent\begin{definition}The edge boundary $b_G(H)$, where $H$
is an induced subgraph of G, is the number of edges which joins
vertices in $H$ to the vertices in $G\backslash H$. Then
$b_G(k)=min\{b_G(H):H\subset G,|H|=k\}$\end{definition}

\begin{theorem}
If $p\geq 2-\sqrt2$, the probability that for some $S\in C_s,
2\leq s \leq 3^{0.44n}$, no edges of $P_{3,p}^n$ join $S$ to
$V(P_3^n)\setminus S$, as $n \rightarrow \infty$, tends to 0.
\end{theorem}
 \proof
 \noindent Consider $S\subset V=V(P_3^n)$ and
set $b(S)=b_{P_3^n}(H)$  where $H$ is the induced subgraph of
$P_3^n$ with the set of vertices S. One can observe that:
\begin{align}
\textbf{P}(\text{No edges of } P_{3,}^n \text{ join S to }
V\setminus S)= (1-p)^{b(S)}
\end{align}
\noindent In order to prove the theorem, it is sufficient to show:
\begin{align}
\sum_{s=2}^{3^{0.44n}}\sum_{S\in C_s}(1-p)^{b(S)}=o(1)
\end{align}

From \cite{edge} one knows that for $|S|=s$:
\begin{align}
b(S)\geq b(s)
\geq\frac{e}{3}s\ln\frac{3^n}{s}=\frac{e\ln3}{3}s(n-\log_3 s)
\end{align}
\noindent one knows that:
\begin{align}
\sum_{S\in C_s}(1-p)^{b(S)}\leq |C_s|(1-p)^{b(s)}
\end{align}

First, one should find a bound for $|C_s|$. One has maximum $3^n$
choices to choose the first element for $C_s$. The selected
element is connected to maximum $2n$ vertices, therefore, there
are maximum $2n$ choices to choose the second element. With the
same logic, there are at most $2n(s-1)$ choices to choose the last
element for $|C_s|$. Hence, one can show:
\begin{align}
|C_s|\leq 3^n(2n)(2n(2))...(2n(s-1)) \leq (s-1)!(2n)^{s-1}3^n
\end{align}
\noindent and:
\begin{align}
|C_s|(1-p)^{b(s)}\leq (s-1)!(2n)^{s-1}3^n
(1-p)^{\frac{e\ln3}{3}s(n-\log_3s)} \label{6.2}
\end{align}
Set $a:=\frac{e\ln3}{3}$. Assume that the right hand sides of
inequality (\ref{6.2}) is equal to A. After multiplying both sides
of inequality (\ref{6.2}) with $\frac{2ns^{s+1}}{s!}$ and then
taking $\log_3$ from both sides, one gets:
\begin{align}
\log_3(|C_s|(1-p)^{b(s)}\frac{2ns^{s+1}}{s!})\leq\log_3(A\frac{2ns^{s+1}}{s!})
\end{align}
\noindent If $\log_3(A\frac{2ns^{s+1}}{s!})\rightarrow-\infty$ as
$n\rightarrow\infty$ then $A\frac{2ns^{s+1}}{s!}$ should tend to 0
for $2\leq s \leq 3^{0.44n}$. This means that
$|C_s|(1-p)^{b(s)}\frac{2ns^{s+1}}{s!}$ tends to 0 for $2\leq s
\leq 3^{0.44n}$, as $n\rightarrow \infty$. Therefore:
\begin{align}
|C_s|(1-p)^{b(s)}\leq \frac{s!}{2ns^{s+1}} \text{ for large values
of n}
\end{align}
which shows that:
\begin{align}
\sum_{s=2}^{3^{0.44n}}\sum_{S\in C_s}(1-p)^{b(S)}=o(1) \label{6.4}
\end{align}

\noindent Finally, it remains to prove
$\log_3(A\frac{2ns^{s+1}}{s!})\rightarrow-\infty$ as
$n\rightarrow\infty$. One can verify this for $s\leq n$ and $s>n$.
After multiplying both sides of inequality (\ref{6.2}) with
$\frac{2ns^{s+1}}{s!}$ and then getting $\log_3$ from both sides,
one gets:
\begin{align}
\log_3(|C_s|(1-p)^{b(s)}\frac{2ns^{s+1}}{s!})\leq
s(\log_32n+\log_3s+an\log_3(1-p)-a(\log_3s)\log_3(1-p))+n
\label{6.3}
\end{align}
For $s\leq n$ the largest factor in equation (\ref{6.3}) is
$n(a(\log_3(1-p))s+1)$ which is negative for $p=p_c$ and $2\leq s
\leq n$. Therefore,
$\log_3(A\frac{2ns^{s+1}}{s!})\rightarrow-\infty$ as
$n\rightarrow\infty$ for $s\leq n$. For $s> n$ the largest factor
in equation (\ref{6.3}) is
$s(an(\log_3(1-p))+(\log_3s)(1-a\log_3(1-p))$ which is negative if
$s<3^{\frac{-a\log_3(1-p)}{1-a\log_3(1-p)}n} \approx 3^{0.44n}$.
Therefore, $\log_3(A\frac{2ns^{s+1}}{s!})\rightarrow-\infty$ as
$n\rightarrow\infty$ for $s>n$.\done\\\\

\noindent \textbf{For $p\geq 0.67$:}\\
\begin{theorem}
If $p\geq 0.67$, the probability that for some $S\in C_s, 2\leq s
\leq 3^{n}$, no edges of $P_{3,p}^n$ join $S$ to
$V(P_3^n)\setminus S$, as $n \rightarrow \infty$, tends to 0.
Hence, the probability that $P^n_{3,p}$ is connected tends to $1$,
as $n \rightarrow \infty$.
\end{theorem}
\proof We know that:
\begin{align}
|C_s|\leq \binom{3^n}{s}<(\frac{e3^n}{s})^s
\end{align}
hence:
\begin{align}
|C_s|(1-p)^{b(s)}\leq
\binom{3^n}{s}<(\frac{e3^n}{s})^s(1-p)^{\frac{e\ln3}{3}s(n-
\log_3s)} \label{6.6}
\end{align}
By using Mathematica \footnote{Mathematica is a computational
software program developed by Wolfram Research of Champaign,
Illinois} one can show the right hand side of equation (\ref{6.6})
tends to zero as $n\rightarrow \infty$.\done

\chapter{Reliable networks}
In many engineering applications, it is of interest to construct a
graph (network), with a specific number of edges and vertices,
which is the most reliable one in that family of graph with n
vertices and m edges. One of the measures of reliability is
all-terminal reliability.
\section{All terminal reliability}
\begin{definition} \textbf{Random subgraph of a graph, Percolation on a
graph:} Random subgraph of $G(V_n,E_m)$ is the graph $G_{p_n}$
which contains all vertices of $G$, and each edge of $G$
independently with probability $p_n$. Doing percolation on a graph
with the parameter $p_n$ is the same as finding a random subgraph
of a graph with the parameter $p_n $. $p_n$ is called percolation
parameter.
\end{definition}
\begin{definition} \textbf{Uniformly optimally reliable graph (UOR)}:
The UOR graph, if it exists, is the graph $G_n\in G(n,m)$ which
maximizes the probability that $G_n$ is connected after
percolation with the parameter $p_n$ for fixed $n,m$ and all
$p_n\in(0,1)$.
\end{definition}
\begin{definition}\textbf{Reliability polynomial}: Let $s_k$
be the number of spanning connected subgraphs of $G_n\in G(n,m)$
having exactly $k$ edges. Let $R(G_n,p_n)$ be the probability that
$G_n$ is connected after percolation with the parameter $p_n$. One
can formulate $R(G_n,p_n)$ as follows:
\begin{align}
R(G_n,p_n):=\sum_{i=0}^m s_i p_n^i(1-p_n)^{m-i} \label{rn.1}
\end{align}
$R(G_n,p_n)$ is called \emph{reliability polynomial} of graph
$G_n$ or \emph{all-terminal reliability}. In this definition,
$s_i=0$ for $i<n-1$ and $s_m=1$. Also, $s_{m-1}$ is $m-$number of
cuts in $G_n$, and $s_{n-1}$ is the number of spanning tress of
$G_n$.
\end{definition}
The UOR graph, if it exists, is the graph $G_n\in G(n,m)$ which
maximizes $R(G_n,p_n)$ for all $p_n\in(0,1)$. From the definition
of reliability polynomial one can see that the value of
reliability polynomial depends on the structure of a graph as well
as percolation value. Trivially, for fixed values of $p_n$ there
always exists an optimal solution for $R(G,p_n)$. In other words,
for fixed values of $p_n$, there always exists a $G_n\in G(n,m)$
which maximizes $R(G_n,p_n)$. The following theorem and corollary
can be helpful to find the UOR graph. This theorem and corollary
are extracted from \cite{survey}.
\begin{theorem} \label{t.1}
Let G and H be two undirected simple graphs both having n nodes
and m edges and $s_k(G),s_k(H)$ denote the number of spanning
connected subgraphs of G and H, respectively, with exactly k edges
\cite{survey}.
\begin{enumerate}
    \item If there exists an integer $0\leq k \leq m-1$ such that
    $s_i(G)=s_i(H)$ for $i=0,1,...,k$ and $s_{k+1}(G)>s_{k+1}(H)$,
    then there exists a $\rho >0$ such that for all $0<p<\rho$ we
    have $R(G,p)>R(H,p)$.
    \item If there exists an integer $0\leq k \leq m$ such that
    $s_i(G)=s_i(H)$ for $i=m,m-1,...,m-k$ and $s_{m-k-1}(G)>s_{m-k-1}(H)$,
    then there exists a $\rho <1$ such that for all $\rho<p<1$ we
    have $R(G,p)>R(H,p)$.
\end{enumerate}
\end{theorem}
\begin{corollary} \label{c.1} If G is UOR, then \cite{survey}:
\begin{enumerate}
    \item G has the maximum number of spanning trees among all
    simple graphs having n nodes and m edges, and
    \item G is $max-\lambda$ , i.e. has the maximum possible value
    of $\lambda$ among all simple graphs having n nodes and m
    edges, namely $\lambda(G)=\lfloor 2m/n \rfloor$, and the
    minimum number of cutsets of size $\lambda$ among all such
    $max-\lambda$ graphs.
\end{enumerate}
where $\lambda(G)$ is the edge connectivity of G, i.e., the
minimum number of edges whose its removal will disconnect $G$.
\end{corollary}
\noindent \textbf{Important coefficients for large n:} Let $n$ be
large and $m$ sufficiently larger than $n$. When $p_n$ is close to
$0$ then $s_{n-1}$, the number of spanning trees of $G_n$, has the
most significant contribution in $R(G_n,p_n)$ since $(1-p_n)$ is
almost $1$ and $p_n^{n-1}$ is much larger than $p_n^m$. Similarly,
when $p_n$ is close to $1$, then $s_{m-1}$, $m-$number of cuts,
has the most significant contribution in $R(G_n,p_n)$.
\subsection*{Laplacian}
The \emph{Laplacian matrix} of a graph is described briefly in
chapter 1. Here, we present a few results on Laplacian. The author
calculated the algebraic connectivity for all graphs with
$n=5,6,7$ and $n-1<m<\binom{n}{2}$  and could not find a direct
relation between algebraic connectivity and all-terminal
reliability. There is a room for further work in this part.

The Laplacian matrix is $L:=(l_{i,j})_{n\times n}$ where
\begin{math}
l:=\left\{%
\begin{array}{ll}
    deg(v_i), & \hbox{if $i=j$;} \\
    -1, & \hbox{if $i\neq j$ and $v_i$ adjacent to $v_j$;} \\
    0, & \hbox{o.w.} \\
\end{array}%
\right.
\end{math}.
Arrange the eigenvalues of $L$ as : $\lambda_1(L) \leq
\lambda_2(L) \leq...\lambda_n(L) $. This set of $\lambda_i$'s are
called the spectrum of $L$ and $\lambda_2(L)$ is called the
algebraic connectivity of a graph. The following interesting lemma
sheds light on some applications of Laplacian matrix (this lemma
is extracted from \cite{Chung}):
\begin{lemma}
For a graph $G$ on $n$ vertices, we have \cite{Chung}\\
\noindent (i): $$\sum_i \lambda_i \leq n$$ with equality holding
if and only if $G$ has no isolated vertices.\\
\noindent (ii):  For $n\geq 2$ $$ \lambda_1 \leq \frac{n}{n-1}$$
with equality holding if and only if $G$ is the complete graph on
$n$ vertices. Also, for a graph $G$ without isolated vertices, we
have: $$ \lambda_{n-1} \geq \frac{n}{n-1}.$$\\
\noindent (iii:) For a graph which is not a complete graph, we
have
$$\lambda_1 \leq 1.$$\\
\noindent (iv:) If $G$ is connected, then $\lambda_1 >0 $. If
$\lambda_i =0$ and $\lambda_{i+1}\neq 0$, then $G$ has exactly
$i+1$ connected components.\\
\noindent (v): For all $i\leq n-1$, we have: $$ \lambda_i \leq 2
$$ with $$\lambda_{n-1}=2$$ if and only of a connected component
of $G$ is bipartite and nontrivial.\\
\noindent (vi:) The spectrum of a graph is the union of the
spectrum of its connected components.
\end{lemma}

\subsection*{Reliability for $n=5,6,7$}
The author calculated $R(G,p)$ for $n=5,6,7$ and
$n-1<m<\binom{n}{2}$. The coefficients of the reliability
polynomials are presented in tables \ref{tab.1}, \ref{tab.2} and
\ref{tab.3}. As table \ref{tab.1} illustrates, for $n=5$ there is
always a UOR graph. For $n=6,7$ there are always a UOR graph
except two cases. If, $(n,m)=(6,11)$, the optimal solution for
$R(G_n,p)$ depends on the value of $p$. For $p<0.29$
(approximately 0.29) graph \ref{fig.1} is optimal while for
$p>0.29$ graph \ref{fig.2} optimal. Also, if $(n,m)=(7,15)$, the
optimal solution for $R(G_n,p)$ depends on the value of $p$. For
$p<0.81$ (approximately 0.81) graph \ref{fig.3} is optimal while
for $p>0.81$ graph \ref{fig.4} is optimal. From these observation
one can conclude that the UOR does not exist for all values of $n$
and $m$.
\begin{table}
\noindent \begin{tabular}{|c|c|c|c|c|c|c|}
  \hline
   n=5 & $s_{m-1}$ & $s_{m-2}$ & ... & ... & ... & $s_{n-1}$ \\
  \hline
  m=5 & 5 &  &  &  &  &  \\
  m=6 & 6 & 12 &  &  &  &  \\
  m=7 & 7 & 20 & 24 &  &  &  \\
  m=8 & 8 & 28 & 52 & 45 &  &  \\
  m=9 & 9 & 36 & 82 & 111 & 75 &  \\
  m=10 & 10 & 45 & 120 & 205 & 222 & 125 \\
  \hline
\end{tabular}\\\\
\caption{\footnotesize The coefficients of the reliability
polynomials of the UOR graph for $n=6$ and $n-1<m<\binom{n}{2}$.}
\label{tab.1}
\end{table}

\begin{table}
\begin{tabular}{|c|c|c|c|c|c|c|c|c|c|c|}
  \hline
  n=6 & $s_{m-1}$ & $s_{m-2}$ & ... & ... & ... & ... & ... & ... & ... & $s_{n-1}$ \\
  \hline
  m=6 & 6 &  &  &  &  &  &  &  &  &  \\
  m=7 & 7 & 16 &  &  &  &  &  &  &  &  \\
  m=8 & 8 & 26 & 36 &  &  &  &  &  &  &  \\
  m=9 & 9& 36 & 78 & 81 &  &  &  &  &  &  \\
  m=10 & 10 & 45 & 116 & 177 & 135 &  &  &  &  &  \\
  \emph{  \textbf{$m_a$=11}} & 11 & 55 & 163 & 309 & 368 & 225 &  &  &  &  \\
\emph{  \textbf{$m_b$=11}} & 11 & 55 & 163 & 310 & 370 & 224 &  &  &  &  \\
  m=12 & 12 & 66 & 220 & 489 & 744 & 740 & 384 &  &  &  \\
  m=13 & 13 & 78 & 286 & 771 & 1249 & 1552 & 1292 & 576 &  &  \\
  m=14 & 14 & 91 & 364 & 999 & 1978 & 2877 & 3040 & 2196 & 864 &  \\
  m=15 & 15 & 105 & 455 & 1365 & 2997 & 4945 & 6165 & 5700 & 3660 & 1296 \\
  \hline
\end{tabular}
\caption{\footnotesize The coefficients of the reliability
polynomials of the UOR graph for $n=6$ and $n-1<m<\binom{n}{2}$.
For $m=11$, and $p<0.29$ (approximately 0.29) the row with $m_a$
is the UOR graph and for $p>0.29$ the row with $m_b$ is the UOR
graph.} \label{tab.2}
\end{table}

\begin{table}
\tiny
\begin{tabular}{|c|c|c|c|c|c|c|c|c|c|c|c|c|c|c|c|}
  \hline
  n=7 & $s_{m-1}$ & $s_{m-2}$ & ... & ... & ... & ... & ... & ... & ... & ... & ... & ... & ... & ... & $s_{n-1}$ \\
  \hline
  m=7 & 7 &   &   &   &   &  &  &  &  &  &  &  &  &  &  \\
  m=8 & 8 & 21 &   &   &  &  &  &  &  &  &  &  & &  &  \\
  m=9 & 9 & 33 & 51 &  &  &  &  &  &  &  &  &  &  & &  \\
  m=10 & 10 & 44 & 104 & 117 &  &  &  &  &  &  &  &  &  &  &  \\
  m=11 & 11 & 55 & 159 & 273 & 231 &  &  &  &  &  &  &  &  &  &  \\
  m=12 & 12 & 66 & 216 & 456 & 612 & 432 &  &  &  &  &  &  &  &  &  \\
  m=13 & 13 & 78 & 284 & 690 & 1146 & 1248 & 720 &   &   &  &  &  &  &  &  \\
  m=14 & 14 & 91 & 364 & 994 & 1932 & 2668 & 2460 & 1200 &  &  &  &  &  &  &  \\
 \emph{ \textbf{$m_a$=15}} & 15 & 105 & 455 & 1360 & 2946 & 4704 & 5464 & 4320 & 1840 &  &  &  &  &  &  \\
\emph{  \textbf{$m_b$=15} }& 15 & 105 & 455 & 1360 & 2946 & 4705 & 5465 & 4305 & 1805 &  &  &  &  &  &  \\
  m=16 & 16 & 120 & 560 & 1817 & 4328 & 7766 & 10548 & 10628 & 7396 & 2800 &   &  &  &  &  \\
  m=17 & 17 & 136 & 680 & 2379 & 6169 & 1226 & 18762 & 22226 & 19808 & 12320 & 4200 &  &  &  &  \\
  m=18 & 18 & 153 & 816 & 3060 & 8562 & 18485 & 31344 & 41964 & 44000 & 35094 & 19716 & 6125 &  &  &  \\
  {\begin{sideways}\parbox{1cm}{\centering m=19}\end{sideways}} &  {\begin{sideways}\parbox{1cm}{\centering 19}\end{sideways}} & {\begin{sideways}\parbox{1cm}{\centering 171}\end{sideways}} & {\begin{sideways}\parbox{1cm}{\centering 969}\end{sideways}} & {\begin{sideways}\parbox{1cm}{\centering 3876}\end{sideways}} & {\begin{sideways}\parbox{1cm}{\centering 11624}\end{sideways}} & {\begin{sideways}\parbox{1cm}{\centering 27073}\end{sideways}} & {\begin{sideways}\parbox{1cm}{\centering 49985}\end{sideways}} & {\begin{sideways}\parbox{1cm}{\centering 73888}\end{sideways}} & {\begin{sideways}\parbox{1cm}{\centering 87468}\end{sideways}} & {\begin{sideways}\parbox{1cm}{\centering 81976}\end{sideways}} & {\begin{sideways}\parbox{1cm}{\centering 58958}\end{sideways}} &  {\begin{sideways}\parbox{1cm}{\centering 30109}\end{sideways}} & {\begin{sideways}\parbox{1cm}{\centering 8575}\end{sideways}} & & \\
  {\begin{sideways}\parbox{1cm}{\centering m=20}\end{sideways}} &  {\begin{sideways}\parbox{1cm}{\centering 20}\end{sideways}} & {\begin{sideways}\parbox{1cm}{\centering 190}\end{sideways}} & {\begin{sideways}\parbox{1cm}{\centering 1140}\end{sideways}} & {\begin{sideways}\parbox{1cm}{\centering 4845}\end{sideways}} & {\begin{sideways}\parbox{1cm}{\centering 15502}\end{sideways}} & {\begin{sideways}\parbox{1cm}{\centering 38725}\end{sideways}} & {\begin{sideways}\parbox{1cm}{\centering 77240}\end{sideways}} & {\begin{sideways}\parbox{1cm}{\centering 124605}\end{sideways}} & {\begin{sideways}\parbox{1cm}{\centering 163400}\end{sideways}} & {\begin{sideways}\parbox{1cm}{\centering 173646}\end{sideways}} & {\begin{sideways}\parbox{1cm}{\centering 147500}\end{sideways}} &  {\begin{sideways}\parbox{1cm}{\centering 96915}\end{sideways}} & {\begin{sideways}\parbox{1cm}{\centering 45530}\end{sideways}} & {\begin{sideways}\parbox{1cm}{\centering 12005}\end{sideways}}& \\
  {\begin{sideways}\parbox{1cm}{\centering m=21}\end{sideways}} &  {\begin{sideways}\parbox{1cm}{\centering 21}\end{sideways}} & {\begin{sideways}\parbox{1cm}{\centering 210}\end{sideways}} & {\begin{sideways}\parbox{1cm}{\centering 1330}\end{sideways}} & {\begin{sideways}\parbox{1cm}{\centering 5985}\end{sideways}} & {\begin{sideways}\parbox{1cm}{\centering 20349}\end{sideways}} & {\begin{sideways}\parbox{1cm}{\centering 54257}\end{sideways}} & {\begin{sideways}\parbox{1cm}{\centering 116175}\end{sideways}} & {\begin{sideways}\parbox{1cm}{\centering 202755}\end{sideways}} & {\begin{sideways}\parbox{1cm}{\centering 290745}\end{sideways}} & {\begin{sideways}\parbox{1cm}{\centering 343140}\end{sideways}} & {\begin{sideways}\parbox{1cm}{\centering 331506}\end{sideways}} & {\begin{sideways}\parbox{1cm}{\centering 258125}\end{sideways}} & {\begin{sideways}\parbox{1cm}{\centering 156555}\end{sideways}} & {\begin{sideways}\parbox{1cm}{\centering 68295}\end{sideways}} & {\begin{sideways}\parbox{1cm}{\centering 16807}\end{sideways}} \\
  \hline
\end{tabular}
\normalsize \caption{\footnotesize The coefficients of the
reliability polynomials of the UOR graph for $n=7$ and
$n-1<m<\binom{n}{2}$. For $m=15$, and $p<0.81$ (approximately
0.81) the row with $m_a$ is the UOR graph and for $p>0.81$ the row
with $m_b$ is the UOR graph.} \label{tab.3}
\end{table}

\begin{figure}
\centering
 \subfigure[\footnotesize ]{
  \includegraphics[scale= 1]{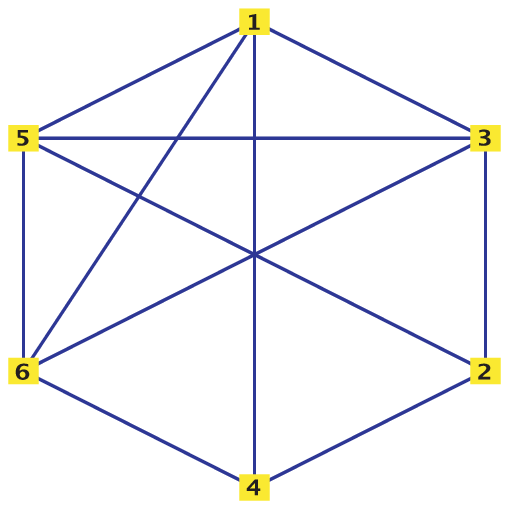}\
 \label{fig.1}
 }\hspace{0.2cm}
\subfigure[\footnotesize ]{
  \includegraphics[scale= 1]{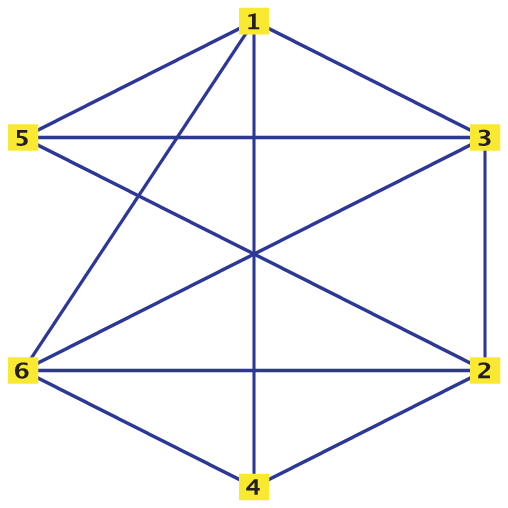}\
  \label{fig.2}
  }
  \caption{\footnotesize $(n,m)=(6,11)$, For $p<0.29$ (approximately 0.29)
graph $(a)$ is optimal while for $p>0.29$ graph $(b)$ is optimal.}
\end{figure}

\begin{figure}
\centering
 \subfigure[\footnotesize ]{
  \includegraphics[scale= 1]{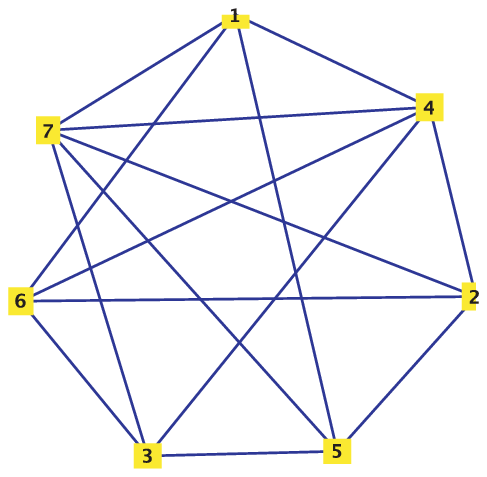}\
 \label{fig.3}
 }\hspace{0.2cm}
\subfigure[\footnotesize ]{
  \includegraphics[scale= 1]{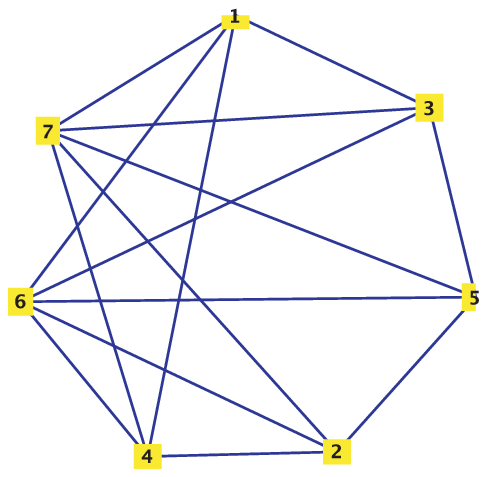}\
  \label{fig.4}
  }
  \caption{\footnotesize $(n,m)=(7,15)$, For $p<0.81$ (approximately 0.81)
graph $(a)$ is optimal while for $p>0.81$ graph $(b)$ is optimal.}
\end{figure}
\subsection*{Reliability for $m=n-1,n,n+1,n+2,n+3$}
For $m=n-1,n,n+1,n+2$ there always exists a UOR graph. For
$m=n-1$, any tree is the UOR graph. For $m=n$, $C_n$, single cycle
with n vertices, is the UOR graph. The first non-trivial case is
$m=n+1$,  which is solved by F. Boesch \cite{Boe1}, \cite{Boe}.
The UOR graph in this case is: for $n\geq 5$, start with a
multigraph with 2 vertices and 3 edges. Then add total of $n-2$
vertices of degree 2 in each lines of the graph so that the number
of vertices in each line differs by at most one \cite{Boe}. For
$m=n+2$ the problem is also solved by F. Boesch. The UOR graph in
this case is: start with $K_4$, then add total of $n-2$ vertices
of degree 2 in each lines of the graph so that the number of
vertices in each line differs by at most one \cite{Boe}. For
$m=n+3$, the UOR graph is found by G. Wang \cite{Wang}. The UOR
graph in this case is: start with $K_{3,3}$, a complete bipartite
graph with 3 vertices in each part, and then add the remanning
vertices as before.
\subsection*{Family of counterexamples}
Kelmans \cite{kel} and Myrvold et al. \cite{Myrvold} found
infinite families of counter examples which the UOR graph does not
exist. As an example, for $n$ even and $n\geq 6$ and
$m=n(n-1)/2-(n+2)/2$, or for $n$ odd and $n\geq 7$ and
$m=n(n-1)/2-(n+5)/2$ there always exists a graph in which it
maximizes $R(G,p_n)$ for $p$ close to $1$, but do not have the
maximum number of spanning trees. Therefore, from theorem
\ref{t.1} and corollary \ref{c.1} the UOR graph does not exists
for these families.

\section{Random accessibility}
 M. Ebneshahrashoob, T. Gao and M. Sobel introduced the concept of random accessibility for simple
 graphs
 \cite{Morteza}. They believe that finding the UOR graph is
 related to the concept of random accessibility. In random
 accessibility, the goal is to find the expectation and the variance
 of the number of transitions $X_j$  needed to visit $j$ new
 vertices in $G\in G(V_n,E_m)$. In this approach the starting
 point is not considered as a new vertex. The result of
 the expectation and the variance can depend on starting point. Hence, one
 should change the starting point depending on the degree of it as a weighing
 factor. If one considers a graph with enough symmetry, the
 result does not depend on starting point.
 From analyzing numerical results, they make the following
 interesting conjuncture:
 \begin{itemize}
    \item If the family of graphs contains both regular and
    non-regular graphs, then the UOR graph is among the regular
    graphs. Also, the expectation for random accessibility of
    graphs is equal or greater than the corresponding result of
    the UOR graph for each value of $j$ close to $m-1$ (with the same ordering of the graph as
    for all-terminal reliability).
 \end{itemize}

\appendix
\renewcommand{\baselinestretch}{1.5}

\end{document}